\documentclass[11pt]{article}
\usepackage[T1]{fontenc}
\usepackage[utf8]{inputenc}
\usepackage{lmodern}
\usepackage[a4paper,margin=28mm]{geometry}
\usepackage{amsmath,amssymb,amsthm,mathtools}
\usepackage{microtype}
\usepackage[hidelinks,hypertexnames=false]{hyperref}
\usepackage{enumitem}
\usepackage[dvipsnames]{xcolor}
\usepackage{graphicx,tikz}
\usetikzlibrary{arrows.meta,calc}

\allowdisplaybreaks[2]
\numberwithin{equation}{section}
\newtheorem{theorem}{Theorem}[section]
\newtheorem{proposition}[theorem]{Proposition}
\newtheorem{lemma}[theorem]{Lemma}

\theoremstyle{definition}

\theoremstyle{remark}
\newtheorem{remark}[theorem]{Remark}
\newcommand{\T}{\mathbb T}
\newcommand{\R}{\mathbb R}
\newcommand{\ip}{\int_{-\pi}^{\pi}}
\newcommand{\io}{\int_{0}^{\pi}}
\newcommand{\Q}{\mathcal Q}

\newcommand{\norm}[1]{\left\|#1\right\|}

\title{Traveling waves and their long-time stability \\ in the Kuramoto Mean Field Game}

\author{Marco Cirant \thanks{Dipartimento di Matematica Tullio Levi-Civita, Universit\`a degli Studi di Padova, Via Trieste 63, 35121 Padova, Italy, {\sf cirant@math.unipd.it}} \quad Mohamed Ghattassi \thanks{Department of Mathematics, New York University Abu Dhabi, Saadiyat Island, P.O. Box 129188, Abu Dhabi, United Arab Emirates, {\sf mg6888@nyu.edu}} \quad Nader Masmoudi \thanks{Department of Mathematics, New York University Abu Dhabi, Saadiyat Island, P.O. Box 129188, Abu Dhabi, United Arab Emirates, Courant Institute School of Mathematics, Computing, and Data Science, New York University, 251 Mercer Street, New York, NY 10012, USA, {\sf nm30@nyu.edu}} \quad Gianmaria Verzini \thanks{Dipartimento di Matematica, Politecnico di Milano, Piazza Leonardo da Vinci 32, 20133 Milano, Italy, {\sf gianmaria.verzini@polimi.it}}}

\date{}

\begin{document}
\maketitle
\begin{abstract}
We study the long-time behavior of the Kuramoto Mean Field Game (MFG), a model of spontaneous 
synchronization among controlled stochastic oscillators. We first show that, above the critical 
interaction strength, the synchronized stationary state is unique up to translations and nondegenerate, and that
there exists a family of small-speed synchronized traveling waves. The presence of these waves is a genuine feature of the game theoretic model, and it is not present in the classical Mean Field one. We then investigate their local stability; in particular, for data close to a reference wave, we prove existence of a solution that stays exponentially close, away from the initial / final conditions, to a wave with possibly different speed, up to a phase correction. The proof combines conservation of momentum, a suitable Lyapunov functional, and a modulation argument. This reveals a new type of stable long-time behavior of solutions in MFG which is not, strictly speaking, of turnpike nature. Indeed, rather than selecting a stationary state, agents decide to collectively travel with (almost) constant speed from initial state towards the one favored by the terminal cost.

\end{abstract}
\noindent

{\footnotesize \textbf{AMS-Subject Classification}}
{\footnotesize 35Q89, 35B40, 35C07, 92B25 
}\\
{\footnotesize \textbf{Keywords}}  
{\footnotesize Mean field games; Kuramoto model; synchronization; traveling waves; long-time behavior; turnpike property}


\section{Introduction}

The main purpose of this paper is to explore the behavior of solutions to the Kuramoto mean field game (MFG) system 
\begin{equation}\label{eq:mfg}
\begin{cases}
 -u_t-\nu u_{xx}+\frac12|u_x|^2= -\kappa\ip\cos(x-y)m(t,y)\,dy , \qquad (t,x)\in(0,T)\times\T ,\\
 m_t-\nu m_{xx}-\partial_x(mu_x)=0,\\
  m(0,\cdot)=m_0,\qquad u(T,\cdot)=u_T,
\end{cases}
\tag{MFG}
\end{equation}
in the long-time regime, with particular focus on transport phenomena. The state $x$ denotes a phase belonging to the 
flat torus $\T=\R/(2\pi\mathbb Z)$, which is represented by $[-\pi,\pi]$, and $\nu > 0$. The system has been 
considered in \cite{CCS23} (see also \cite{YMMS12}) to describe Nash equilibria among a population of homogeneous stochastic, controlled oscillators in mean field interaction. It has been inspired by the classical Kuramoto model (see for instance \cite{strogatz}), which is a paradigmatic model of spontaneous synchronization. Its key feature is a phase transition from incoherent states to self-organization when the interaction parameter $\kappa > 0$ exceeds a certain threshold. 

The classical stochastic Kuramoto model has been extensively studied in the mean field limit, and its long-time behavior is 
now well understood \cite{GPP12}: using its gradient flow structure, it is known that for all values of $\kappa$, solutions 
must converge to a stationary state as $t \to \infty$; here the characteristic frequency $\omega$ of the oscillators is indeed normalized to zero. Although for the game-theoretic analog \eqref{eq:mfg} 
considered here very similar phase transition phenomena  have been observed as the parameter $\kappa$ varies, the study of its long-time behavior is still largely open. One main difficulty is that \eqref{eq:mfg} is a forward-backward system, and uniqueness of solutions is false in general. The question we address in this work is then the following: given initial and terminal conditions $m_0, u_T$, does there exist 
at least one solution that exhibits specific long-time patterns? Under monotonicity conditions, MFGs are known to converge to 
stationary states as $T \to \infty$. However, the model considered here is not monotone and there is no a priori reason why every dynamic  equilibrium should become stationary in the long run.

\medskip 
The emergence of traveling waves and periodic solutions has long attracted interest in game theory, particularly in MFG. Already in the lectures of Lions \cite{LionsSeminar}, situations in which equilibria exhibit this kind of pattern, such as the “Mexican wave”, were discussed. These patterns have been observed in various models that incorporate some form of aggregation or imitation among players \cite{CC21, C19, PR22}. Traveling waves have also previously been observed in the Kuramoto MFG model itself \cite{YMMS12}. The aim of this work is to determine whether these kind of equilibria are robust: if, for instance, the initial conditions are perturbed, can the players still reach an equilibrium exhibiting the same pattern? Surprisingly, traveling waves are absent in the analogous classical Kuramoto model with identical oscillators, which is forward in time. Their presence here indicates  spontaneous synchronization at frequencies that are \emph{different} from the natural one (which is zero).

\medskip

As for the standard model, the game-theoretic one has several applications, see for example \cite{MM25}. In~\cite{MR4064663}, oscillators are neurons that adjust their phases in response to neighboring cells and environmental cues, balancing the cost of modifying their natural frequencies against alignment. In \cite{Getal26}, a Kuramoto MFG model is used to study how incentives can affect synchronization of an electrical grid; each agent chooses how much power to inject or withdraw under
a frequency-dependent price. In both situations, traveling waves constitute collective periodic motions that arise from purely individual decisions.



\medskip

Let us briefly summarize first what is known on the long-time behavior of MFG, and the MFG Kuramoto model. As anticipated above, whenever \eqref{eq:mfg} possesses a monotone structure, it is now well-known that dynamic and stationary solutions are unique, and the former converge to the latter typically exponentially fast, except close to $t = 0$ and $t = T$, where initial and final conditions are imposed (as in \eqref{eq:intropike} below). This behavior is typically referred to as the \emph{turnpike property}, which is a broadly studied phenomenon in optimal control \cite{TZSur}. See \cite{CLLP1, CM24, CP21} and references therein for related results in MFG. The Kuramoto MFG falls out of any known monotone regime, and uniqueness of solutions is actually false for $\kappa$ large; this is clearly seen by looking at the stationary problem, where any non-trivial solution generates a continuum of solutions in view of the invariance by $x$-translation of the system. 

The first dynamic analysis of the Kuramoto MFG has been carried out in \cite{CCS23}, in the discounted case; here the authors proved that the  transition to synchronization occurs in the MFG framework as a bifurcation from the \emph{incoherent state} $\bar m = (2\pi)^{-1}$ at a critical threshold $\kappa = \kappa_c$, with full synchronization occurring in the $\kappa \to \infty$ limit. Moreover, in the sub-critical range $\kappa < \kappa_c$, the incoherent state was shown to be locally stable in long time. This has been further generalized to heterogeneous models with intrinsic frequencies \cite{CCS25}. The global stability of the incoherent state in the subcritical regime has been then obtained in \cite{hofer2026convergence}. The recent paper \cite{munoz2026uniqueness} has provided a complete picture of the \emph{stationary} solutions of the discounted problem: while up to the critical threshold the only solution is the incoherent one, above it synchronized solutions are unique up to phase translation (see Theorem~\ref{thm:classification} below for the stationary solutions of \eqref{eq:mfg}).

The first result on the long time stability of \emph{synchronized states} has been obtained in \cite{cesaroni2024stationary}. Here, it was shown that for large enough $\kappa$, an even synchronized state is uniquely determined, and it locally attracts \emph{even} solutions of the dynamic problem \eqref{eq:mfg}. In a sense, restricting to the symmetric setting restores a sort of (local) uniqueness of stationary solutions. As global uniqueness of dynamic solutions cannot be expected, a local turnpike property of the following kind was established:  denoting by $(U,M)$ the unique (even) synchronized state, if $(u(t),m(t))$ is close enough to $(U,M)$ and even, then it must satisfy
\begin{equation}\label{eq:intropike}
\|m(t)-M\| + \|u(t)-U\| \lesssim e^{-\omega t} + e^{-\omega (T-t)}
\end{equation}
for some $\omega$ not depending on $T$.

One of the main goals of this work is to establish analogous long-time convergence results without imposing any symmetry assumptions on the data. In this more general setting, the landscape of stationary states becomes richer, consisting of a continuum of translates of the synchronized state. Moreover, the presence of families of \emph{traveling waves} comes into the scene; as we shall see later, their thorough study will be crucial to understanding the long-time behavior of the system in the absence of symmetry assumptions on the data.

\medskip

Let us state explicitly our results on traveling waves, introducing first some notations. For $\rho\in L^1(\T)$, we define the nonlocal coupling \[
 K[\rho](x)=-\kappa\ip\cos(x-y)\rho(y)\,dy = -\kappa\bigl(C[\rho]\cos x+S[\rho]\sin x\bigr),
\]
and for brevity,
\begin{equation}\label{eq:CSR}
 C[\rho]=\ip\rho\cos x\,dx,\qquad
 S[\rho]=\ip\rho\sin x\,dx,\qquad
 R[\rho]=(C[\rho]^2+S[\rho]^2)^{1/2}.
\end{equation}
Traveling waves will be solutions of the form
\begin{equation}\label{eq:def_TW}
 u(t,x)=U(x-\xi t)+\bar\lambda(T-t),\qquad m(t,x)=M(x-\xi t),
\end{equation}
where $\xi\in\R$ is the speed of the wave (the case $\xi=0$ corresponds to
stationary solutions) and the triple $(U,M,\bar\lambda)$ solves a suitable 
stationary system, see \eqref{eq:profiles} below. For every $\xi$ there is the \emph{incoherent}
solution, with $U\equiv0$ and $M\equiv(2\pi)^{-1}$, while a solution with
$R[M]>0$ is called \emph{synchronized}, and it gives rise to a nontrivial wave (we refer to
Section~\ref{sec:stationary} for more details). We have the following
classification.

\begin{theorem}\label{thm:classification}
Let $\kappa_c(\xi)=2(\nu^2+\xi^2)$, for $\xi\in\R$, and
$\kappa_c=\kappa_c(0)=2\nu^2$. Then:
\begin{enumerate}
\item[(a)] Let $\xi=0$. If $\kappa\le\kappa_c$, the incoherent solution is the
 only stationary solution. If $\kappa>\kappa_c$, there is a synchronized
 stationary solution $(U_0,M_0,\bar\lambda_0)$, which is unique up to
 translations and nondegenerate, in the sense of \eqref{eq:radial-gap} below.
 Moreover, the map $\kappa\mapsto R_0(\kappa):=R[M_0]$ is smooth and strictly increasing on
 $(\kappa_c,\infty)$, and $R_0(\kappa)\to0$ as $\kappa\downarrow\kappa_c$.
\item[(b)] For every $\kappa>0$ there exists $\beta_0>0$, depending only on
 $\kappa/\nu^2$, such that the following holds for every $|\xi|\le\nu\beta_0$.
 If $\kappa\le\kappa_c(\xi)$, the incoherent solution is the only traveling
 wave of speed $\xi$. If $\kappa>\kappa_c(\xi)$, there is a synchronized
 traveling wave of speed $\xi$, which is unique up to translations. We denote by $(U_\xi,M_\xi,\bar\lambda_\xi)$ such a wave, normalized by $C[M_\xi]>0$, $S[M_\xi]=0$.
\end{enumerate}
\end{theorem}

The classification result in point (a) is a consequence of concavity of a certain fixed-point map 
(sometimes called ``response map'', or ``self-consistency map'' in \cite{munoz2026uniqueness}). Such a 
concavity property, first conjectured in \cite[Rem.~7.4]{CCS23}, has been proved recently by Munoz 
\cite{munoz2026uniqueness} in the discounted case. He also anticipates that his method should apply to 
the ergodic problem. Our proof of
Proposition~\ref{prop:concavity} below follows his strategy, with some simplifications. In the ergodic 
case the second derivative of the fixed-point map, see \eqref{eq:A2} below, contains only one of his two 
terms, thus the ``cubic-moment'' part of \cite{munoz2026uniqueness} is not needed. The concavity 
property can be also seen as the convexity of the derivative of the lowest characteristic value of the 
\emph{Mathieu equation}, which is a classical and widely studied ODE in the literature. We record this 
observation in Remark \ref{rem:mathieu}.

A further observation here is the nondegeneracy property. This allows, by a perturbation argument, to 
construct genuine traveling waves, that is point (b) of the previous theorem. Furthermore, it is 
fundamental to obtain that a certain quadratic form is positive definite, orthogonal to the translation 
direction, which will be crucial in the analysis of the time-dependent problem.

\medskip

We now turn to the main long-time convergence result, whose statement is as follows (it is basically Theorem \ref{thm:continuation} below).

\begin{theorem}\label{thm:existenceintro} Let $\kappa>\kappa_c$, and let $(U_\xi,M_\xi,\bar\lambda_\xi)$ denote the synchronized traveling waves of Theorem~\ref{thm:classification}. Assume that $m_0,u_T\in C^{2,\alpha}(\T)$ for some $\alpha \in (0,1)$, $m_0>0$, $\ip m_0\,dx=1$. There exist $e_0, \xi_0, T_0 > 0$ (depending on $\nu, \kappa$) such that, whenever $T \ge T_0$ and
\[
\|m_0-M_{\xi}(\cdot-a)\|_{L^2} +\|\partial_x u_T-U'_{\xi}(\cdot-\xi T-a)\|_{L^2} \le e_0
\]
for some $a$ and $|\xi| \le \xi_0$, then \eqref{eq:mfg} has at least one classical solution $(u,m)$ satisfying the following property : there exists $\xi_*$ and $q \in C^1([0,T])$ such that
\[
|\xi_*-\xi| \le \frac{C}{T}e_0
\]
and
\begin{multline}\label{eq:expintro}
|q'(t)|+\|m(t, \cdot)-M_{\xi_*}(\cdot-\xi_* t -q(t))\|_{L^2}+\|u_x(t, \cdot)-U'_{\xi_*}(\cdot-\xi_* t -q(t))\|_{L^2} \le \\ Ce_0\left(e^{-\omega t}+e^{-\omega(T-t)}\right) \qquad \forall t \in [0,T]
\end{multline}
for some $C, \omega > 0$ independent of $T$.
\end{theorem}

Roughly speaking the theorem says that, if initial and final conditions are close to a reference wave 
with speed $\xi$, then there exists a solution of \eqref{eq:mfg} that is exponentially close to 
another wave, with speed $\xi_*$ not necessarily equal to $\xi$, and subject to a further phase 
correction $q(t)$ (whose speed is exponentially small). This is a sort of stability result of the waves 
obtained in Theorem \ref{thm:classification}, that may be reminiscent 
of stability phenomena that have been observed in several finite and infinite 
dimensional Hamiltonian systems \cite{MR1104946,MR338584,GSS87} (see also 
\cite{PegoWeinstein94}, where the limiting wave has a slightly different speed and phase, and 
\cite{Sattinger76} for the stability with asymptotic phase of traveling waves in parabolic systems). It is not, strictly speaking, of turnpike nature, as usually happens in MFG.

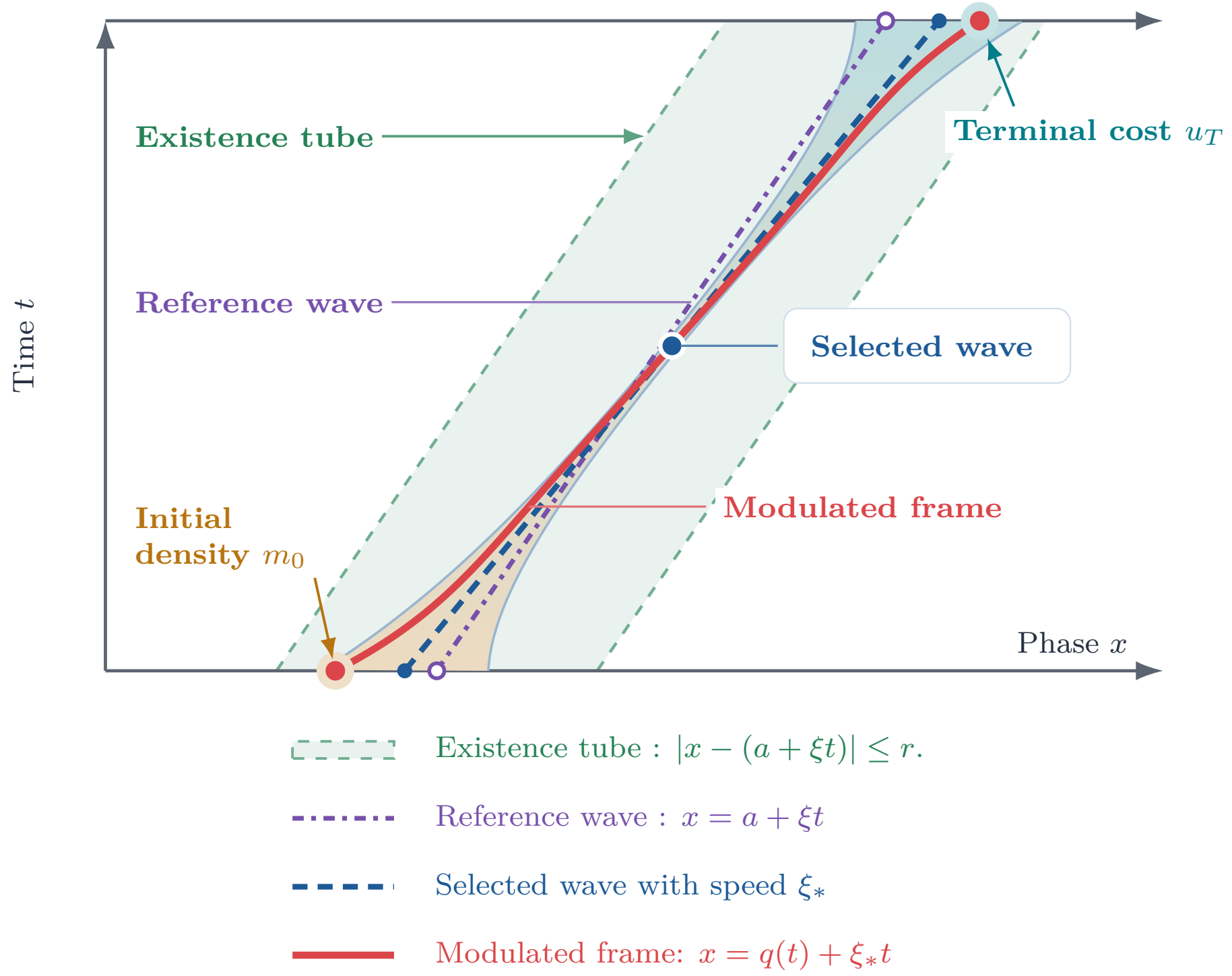
\begin{figure}[!htbp]
\centering
\begingroup
\definecolor{ktpblue}{RGB}{30,90,151}
\definecolor{ktpcoral}{RGB}{219,69,73}
\definecolor{ktpviolet}{RGB}{118,80,171}
\definecolor{ktpamber}{RGB}{183,115,15}
\definecolor{ktpteal}{RGB}{0,125,137}
\definecolor{ktpgreen}{RGB}{37,130,87}
\definecolor{ktpink}{RGB}{38,48,65}

\resizebox{!}{.40\textheight}{%
\begin{tikzpicture}[
  x=1cm,y=.80cm,font=\footnotesize,
  >=Latex,
  axis/.style={
    draw=ktpink!75,line width=.85pt,
    -{Latex[length=2.6mm,width=1.8mm]}
  },
  selected/.style={
    draw=ktpblue,line width=1.5pt,
    dash pattern=on 4pt off 2.7pt
  },
  phase/.style={
    draw=ktpcoral,line width=1.9pt
  },
  prescribed/.style={
    draw=ktpviolet,line width=1.25pt,dash dot
  }
]

\fill[ktpgreen!10]
  (0,1.6)--(4.2,9.2)--(7.2,9.2)--(3,1.6)--cycle;
\draw[ktpgreen!65,dashed,line width=.8pt]
  (0,1.6)--(4.2,9.2);
\draw[ktpgreen!65,dashed,line width=.8pt]
  (3,1.6)--(7.2,9.2);

\draw[axis] (-1.6,1.6)--(-1.6,9.2);
\node[rotate=90,text=ktpink,font=\scriptsize]
  at (-2.37,5.4) {Time $t$};


\draw[axis] (-1.6,1.6)--(8.30,1.6);
\draw[axis] (-1.6,9.2)--(8.30,9.2);

\node[anchor=east,text=ktpink,font=\scriptsize]
  at (8.10,1.93) {Phase $x$};


\begin{scope}
\clip
  plot[domain=0:1,samples=120]
    ({1.2+5.0*\x+.78*(exp(-6*\x)+exp(-6*(1-\x)))},
     {1.6+7.6*\x})
  -- plot[domain=1:0,samples=120]
    ({1.2+5.0*\x-.78*(exp(-6*\x)+exp(-6*(1-\x)))},
     {1.6+7.6*\x})
  -- cycle;
\shade[
  bottom color=ktpamber!26,
  middle color=ktpblue!10,
  top color=ktpteal!26
]
  (0,1.6) rectangle (8.3,9.2);
\end{scope}

\draw[ktpblue!45,line width=.65pt]
  plot[domain=0:1,samples=120]
    ({1.2+5.0*\x+.78*(exp(-6*\x)+exp(-6*(1-\x)))},
     {1.6+7.6*\x});
\draw[ktpblue!45,line width=.65pt]
  plot[domain=0:1,samples=120]
    ({1.2+5.0*\x-.78*(exp(-6*\x)+exp(-6*(1-\x)))},
     {1.6+7.6*\x});

\draw[prescribed] (1.5,1.6)--(5.7,9.2);
\draw[selected] (1.2,1.6)--(6.2,9.2);
\draw[phase]
  plot[domain=0:1,samples=180]
    ({1.2+5.0*\x
      -.38*(exp(-6*\x)-exp(-6*(1-\x)))
      -.27*exp(-6*\x)*(1-\x)*(1-2*\x)
      +.16*sin(720*\x)*(exp(-6*\x)+exp(-6*(1-\x)))},
     {1.6+7.6*\x});

\fill[ktpamber!22] (.5509,1.6) circle (5pt);
\fill[ktpteal!22] (6.5791,9.2) circle (5pt);
\fill[ktpcoral] (.5509,1.6) circle (2.6pt);
\fill[ktpcoral] (6.5791,9.2) circle (2.6pt);
\fill[ktpblue] (1.2,1.6) circle (2pt);
\fill[ktpblue] (6.2,9.2) circle (2pt);
\draw[ktpviolet,fill=white,line width=1pt]
  (1.5,1.6) circle (2pt);
\draw[ktpviolet,fill=white,line width=1pt]
  (5.7,9.2) circle (2pt);

\fill[white] (3.70,5.4) circle (3.7pt);
\fill[ktpblue] (3.70,5.4) circle (2.5pt);

\node[
  anchor=west,align=left,font=\scriptsize,
  text=ktpblue,fill=white,draw=ktpblue!20,
  rounded corners=1.3mm,inner sep=2.5mm
] at (4.74,5.4) {
  \textbf{Selected wave}
};
\draw[ktpblue!75,line width=.6pt]
  (4.69,5.4)--(3.78,5.4);

\node[
  anchor=west,text=ktpviolet,font=\scriptsize,
  fill=white,inner sep=2pt
] (prescribedlabel) at (-1.40,5.91) {
  \textbf{Reference wave}
};
\draw[ktpviolet!75,line width=.6pt]
  (prescribedlabel.east)--(3.88,5.91);

\node[
  anchor=west,align=left,text=ktpcoral,
  font=\scriptsize,fill=white,inner sep=2pt
] at (4.10,3.52) {
  \textbf{Modulated frame}
};
\draw[ktpcoral!75,line width=.6pt]
  (4.02,3.52)--(2.3608,3.52);

\node[
  anchor=west,text=ktpgreen,font=\scriptsize,
  fill=white,inner sep=2pt
] (existencelabel) at (-1.40,7.85) {
  \textbf{Existence tube}
};
\draw[ktpgreen!75,-{Latex[length=2mm]},line width=.7pt]
  (existencelabel.east)--(3.45,7.85);

\node[
  anchor=west,align=left,text=ktpamber,
  font=\scriptsize,fill=white,inner sep=2pt
] (initialdensity) at (-1.40,3.15) {
  \textbf{Initial}\\ \textbf{density $m_0$}
};
\draw[ktpamber,-{Latex[length=2mm]},line width=.7pt,shorten >=3pt]
  (initialdensity.south east)--(.5509,1.6);

\node[
  anchor=west,align=left,text=ktpteal,
  font=\scriptsize,fill=white,inner sep=2pt
] at (6.25,7.90) {
  \textbf{Terminal cost $u_T$}
};
\draw[ktpteal,-{Latex[length=2mm]},line width=.7pt]
  (6.9,8.2)--(6.65,9);

\begin{scope}[yshift=5mm]
\filldraw[
  fill=ktpgreen!10,draw=ktpgreen!65,
  dashed,line width=.7pt
]
  (.15,-.05) rectangle (1.10,.15);

\node[
  anchor=west,align=left,
  font=\scriptsize,text=ktpgreen
] at (1.35,.05) {
  Existence tube : $|x-(a+\xi t)|\le r$.
};

\draw[prescribed] (.15,-.75)--(1.10,-.75);
\node[anchor=west,font=\scriptsize,text=ktpviolet]
  at (1.35,-.75) {
    Reference wave : $x=a+\xi t$
  };

\draw[selected] (.15,-1.55)--(1.10,-1.55);
\node[anchor=west,font=\scriptsize,text=ktpblue]
  at (1.35,-1.55) {
    Selected wave with speed $\xi_*$
  };

\draw[phase] (.15,-2.35)--(1.10,-2.35);
\node[anchor=west,font=\scriptsize,text=ktpcoral]
  at (1.35,-2.35) {
    Modulated frame: $x=q(t)+ \xi_* t$
  };

\end{scope}

\end{tikzpicture}%
}

\caption{Schematic overview of Theorem~\ref{thm:existenceintro}.}\label{figure}
\label{figmainturnpikesummary}
\endgroup
\end{figure}

The new phenomenon emerging here is already important for the understanding of the long-term behavior
of the even synchronized stationary state $(U_0,M_0)$. 
Assume that $(m_0, u_T)$ is a small perturbation of $(U_0,M_0)$, not necessarily even. One may expect 
that, in long time, $(u(t), m(t))$ is exponentially close to $(U_0,M_0)$, possibly up to translation, in 
line with standard turnpike results. By contrast, what happens here is that $(u(t), m(t))$ \emph{is transported} with speed $|\xi_*| \lesssim 1/T$, in order to connect the initial condition to the final 
one, which is favored by the terminal cost. A sort of turnpike behavior is revealed by selecting the 
appropriate moving frame. It should be noticed that, if $(m_0, u_T)$ also is assumed to be even, we 
recover the exponential turnpike property close to $(U_0,M_0)$, as discussed in Remark 
\ref{rem:momentum} below. Moreover, the selected speed $\xi_*$ is uniquely determined by the argument 
described below, but the exponential estimate \eqref{eq:expintro} is robust under small perturbations of 
$\xi_*$ and $q$, see Remark \ref{rem:uniquexi}.
A schematic picture of the result is shown in Figure \ref{figure}.


\medskip

The proof of the existence theorem is given in Section \ref{sec:tdp}, and uses a continuation 
argument that exploits an a priori estimate on solutions $(u,m)$ that are close to the reference 
wave $(U_{\xi},M_{\xi})$. The first crucial observation is that the momentum $\ip mu_x\,dx$ is constant in time, but may differ from the one of the reference wave. It is then natural to compare $(u,m)$ with another wave $(U_{\xi_*},M_{\xi_*})$ that matches the momentum of 
$(u,m)$. When speeds are small, this is uniquely determined by monotonicity of the map $\xi \mapsto \int 
M_\xi U'_\xi$.

Then, we look at the evolution of the quantity
\[
\mathcal L(t) = \ip(u(t)-U_{\xi_*})(m(t)-M_{\xi_*})\,dx + \nu \ip\frac{(m(t)-M_{\xi_*})^2}{M_{\xi_*}}\,dx
\]
This is the standard state-costate product plus a correction term that should be considered as the second 
variation of the relative entropy of $m$ with respect to $M_{\xi_*}$. Looking at this quantity turns 
out to be effective in the study of the long-time behavior of \emph{quadratic} MFG \cite{CD26, 
hofer2026convergence}. In particular, arguing as in \cite{CD26}, one expects that
\begin{equation}\label{eq:introL}
\mathcal L(t)' \le -\mathcal Q(t) + \text{t.o.t.}
\end{equation}
where $\mathcal Q(t)$ is a quadratic form associated to the stationary problem, and ``t.o.t.'' are 
negligible third order terms. A major issue here is that $\mathcal Q(t)$ is not positive definite, as  stationary solutions are translation invariant. The analysis of the first part of the work then comes 
into play, and shows that, orthogonal to this direction of degeneracy, $\mathcal Q(t)$ is positive 
definite. Then a \emph{modulation argument} is developed, that is, $(U_{\xi_*},M_{\xi_*})$ is 
evaluated at $x - t \xi_* - q(t)$, where $q$ is a further phase correction, in order to have $(u(t)-
U_{\xi_*}, m-M_{\xi_*})$ orthogonal to the bad direction. While this restores coercivity of $\mathcal 
Q$, it introduces additional quadratic terms in \eqref{eq:introL}. The last step is then to show that 
those terms can be made small, and therefore can be reabsorbed by $\mathcal Q$ itself, whenever $\xi_*$ 
is small enough; here, the initial choice of $\xi_*$ matching exactly the momentum of the solution with 
the one of the wave is crucial.

Having finally $\mathcal L$ decreasing in time, hence playing the role of a Lyapunov function, gives 
the exponential convergence of $(u,m)$ to $(U_{\xi_*},M_{\xi_*})$. A further technical step shows that  while $(u,m)$ is initially supposed to be close to the reference wave $(U_{\xi},M_{\xi})$, 
it belongs in fact to a smaller neighborhood, and this allows to set up a topological continuation 
argument to get existence.

\medskip

Some final remarks are now in order. First of all, by periodicity in the $x$-variable, any 
traveling wave can be thought of as a \emph{periodic in time} solution.  To our knowledge, results 
obtained here are the first ones showing that those solutions are stable in some sense, and that they describe the long time behavior even in situations where there is no explicit transport effect. It is indeed agents' decisions that create 
transport to connect misaligned initial / final conditions. These equilibria survive as $T \to \infty$, giving rise to families of periodic solutions with arbitrary speed. Further discussion on this aspect is given in Remark 
\ref{rem:periodic}.

Second, note that all the results here involve traveling waves with a small speed, and there are several  places where smallness comes into play. For instance, their uniqueness and nondegeneracy are obtained by analogous properties of the stationary problem, and these follow in turn by concavity of the fixed-point map. Such a map fails to be concave if the speed is large enough, and genuine multiplicity of waves arises, as discussed in Remark \ref{rem:kappac}. Moreover, small speed is required in the dynamic stability argument, as described above. Hence, the stability analysis 
of the synchronized states in the more general model with intrinsic frequencies \cite{CCS25}, where 
multiple waves with arbitrarily large speed may interact, presents a challenging problem at this stage.

We would also like to point out that the analysis is developed here for a specific model, but the 
techniques are flexible enough to cover other MFG problems with the following features : lack of  monotonicity, translation invariance, and conservation of certain quantities (such as the momentum). An 
important structural assumption is the quadratic nature of the Hamiltonian; the Lyapunov function used 
here might not be suitable for non quadratic problems, see the discussion in \cite{CD26}. The analysis of the long-time behavior of MFG beyond classical monotone frameworks has recently received considerable attention \cite{AM26, BK, GKM26, LMN26, PR22}, and it would be interesting to investigate the dynamics of the traveling wave type in different models.

\medskip

Finally, we mention that the modulation method, developed here in the context of MFG, has been extensively 
used to prove stability properties and to characterize the long-time asymptotic dynamics of various nonlinear PDEs. For 
example, in \cite{MR3861072} modulation is used to determine how fast a density of critical mass concentrates in the 
two-dimensional Patlak-Keller-Segel equation. For the steady Prandtl equation, in \cite{MR4028516} the wall shear is used 
as a modulation parameter to study separation under an adverse pressure gradient; see also \cite{MR4493613} for results 
on the unsteady Prandtl flow.

\medskip

The paper is organized as follows. Section~\ref{sec:stationary} 
is devoted to stationary solutions and traveling waves: after the reduction to a scalar fixed point 
problem, we prove Theorem~\ref{thm:classification}, relying on the concavity of the fixed-point map, and 
we establish the coercivity properties of the waves which are needed in the dynamical analysis. In 
Section~\ref{sec:tdp} we deal with the time-dependent problem: we first obtain the exponential estimates 
for the solutions which are close to a reference wave (Theorem~\ref{thm:half}), then we show that the 
linearized system has only the trivial solution, and finally we prove Theorem~\ref{thm:existenceintro} 
by a continuation argument (Theorem~\ref{thm:continuation}). Appendix~\ref{sec:concavity} contains the 
proof of the concavity of the fixed-point map in the stationary case.

\bigskip
	
	\textbf{Acknowledgements.} M. C. and G. V. are members of the Gruppo Nazionale per l'Analisi Matematica, la Probabilit\`a e le loro Applicazioni (GNAMPA) of the Istituto Nazionale di Alta Matematica (INdAM). M. C. has been partially funded by the King Abdullah University of Science and Technology Research Funding (KRF) under award no. CRG2024-6430.6. 
	
\bigskip
	
	\textbf{AI use declaration.} AI tools were used to proofreading the manuscript and to assist in simplifying and refining the presentation of some proofs. The authors assume full responsibility for the mathematical content and the references.

\section{Stationary solutions, traveling waves, and their coercivity}
\label{sec:stationary}

This section is devoted to the stationary solutions of \eqref{eq:mfg} and to
its traveling waves of small speed, which are treated in a unified way
(Theorem~\ref{thm:waves}). The properties of the traveling waves used in the
rest of the paper are collected in Theorem~\ref{thm:travel-gap}, while the
proof of Theorem~\ref{thm:classification} is given in Section~\ref{subsec:concavity}.

Taking into account \eqref{eq:def_TW}, we call a  \emph{traveling wave} of \eqref{eq:mfg}, with speed $\xi\in\R$, 
any triple $(U,M,\bar\lambda)$ satisfying
\begin{equation}\label{eq:profiles}
\begin{cases}
 -\nu U''+\tfrac12|U'|^2+\xi U'+\bar\lambda=K[M],\\
 -\nu M''-\bigl((U'+\xi)M\bigr)'=0,\\
 M>0,\qquad \ip M\,dx=1,\qquad U(0)=0 ;
\end{cases}
\end{equation}
in case $\xi=0$, it is a \emph{stationary solution}. In particular, the
associated \emph{flux} $J:=\nu M'+(U'+\xi)M$ is constant and, dividing by $M$
and integrating on $\T$, we obtain
\begin{equation}\label{eq:flux}
 \nu M'+(U'+\xi)M=J,\qquad\text{where}\quad J\ip M^{-1}\,dx=2\pi\xi ;
\end{equation}
thus $J=0$ if and only if $\xi=0$.

The triple $(0,(2\pi)^{-1},0)$ solves \eqref{eq:profiles} for every
$\xi\in\R$: it is the \emph{incoherent} solution; a solution
with $R[M]>0$ is called \emph{synchronized} (recall that the integral coefficient $R[\rho]$ is defined in 
\eqref{eq:CSR}, together with $C[\rho]$, $S[\rho]$). We immediately notice that,
if $(U,M,\bar\lambda)$ solves \eqref{eq:profiles}, then for some $q\in\T$,
$C[M]=R[M]\cos q$, $S[M]=R[M]\sin q$, so that 
\[
 K[M](x)=-\kappa R[M]\cos(x-q),\qquad
 C[M(\cdot+q)]=R[M],\qquad S[M(\cdot+q)]=0,
\]
and $(U(\cdot+q)-U(q),M(\cdot+q),\bar\lambda)$ is still a solution. Then, up
to translations, we can always reduce to the case $C[M]=R[M]$, $S[M]=0$. 
Finally, for $\xi\in\R$ we write
\begin{equation}\label{eq:kappac}
 \kappa_c(\xi)=2(\nu^2+\xi^2),\qquad \kappa_c=\kappa_c(0)=2\nu^2 ,
\end{equation}
as in Theorem~\ref{thm:classification} (see also Remark~\ref{rem:kappac} below).

The main goal of this section is to show the following result, concerning
the synchronized traveling waves of small speed (stationary solutions included).

\begin{theorem}\label{thm:waves}
Let $\kappa>\kappa_c$. There exists $\beta_0>0$ such that the following holds
for every speed $\xi\in\R$ satisfying
\begin{equation}\label{eq:speeds}
 |\xi|\le\nu\beta_0\qquad\text{and}\qquad \kappa>\kappa_c(\xi).
\end{equation}
\begin{enumerate}
\item[(i)] There is a synchronized solution
 $(U_\xi,M_\xi,\bar\lambda_\xi)$ of \eqref{eq:profiles} with 
 \[
 C[M_\xi]>0,\qquad S[M_\xi]=0,
 \] 
and every solution of \eqref{eq:profiles} is either incoherent or a translate
 \[
  \bigl(U_\xi(\cdot-q)-U_\xi(-q),\ M_\xi(\cdot-q),\ \bar\lambda_\xi\bigr),
  \qquad q\in\T .
 \]
\item[(ii)] The map $\xi\mapsto(U_\xi,M_\xi,\bar\lambda_\xi)$ is smooth on the
 (open) set of the speeds satisfying \eqref{eq:speeds}, with values in
 $C^k(\T)^2\times\R$, for every $k$.
\item[(iii)] $M_\xi$ is even, and
 \[
  M_{-\xi}=M_\xi,\qquad U_{-\xi}(x)=U_\xi(-x),\qquad
  \bar\lambda_{-\xi}=\bar\lambda_\xi .
 \]
\item[(iv)] For $\xi=0$, let $R_0=C[M_0]$. Then the stationary solution
 $(U_0,M_0,\bar\lambda_0)$ is nondegenerate, in the sense that
 \begin{equation}\label{eq:radial-gap}
  \delta_0:=1-\frac{\kappa}{\nu^2}\ip M_0|g_0'|^2\,dx>0,
 \end{equation}
 where $g_0\in H^1(\T)$ is the solution of
 $-(M_0g_0')'=M_0(R_0-\cos x)$ with $\ip M_0g_0\,dx=0$.
\end{enumerate}
\end{theorem}

\subsection{Reduction to a scalar problem}\label{subsec:reduction}

The proof of Theorem \ref{thm:waves} is based on the well known reduction of 
\eqref{eq:profiles} to a fixed point problem. To obtain it, we first study the following 
auxiliary problem. For $\tau\ge0$ and $\beta\in\R$ consider the problem
\begin{equation}\label{eq:scalar}
\begin{cases}
 -v''+\tfrac12|v'|^2+\beta v'+\Lambda=\tau(1-\cos x),\qquad v(0)=0,\\[2pt]
 -m''-\bigl((v'+\beta)m\bigr)'=0,\qquad m>0,\quad \ip m\,dx=1 .
\end{cases}
\end{equation}

Solutions of \eqref{eq:scalar} will be denoted with 
$(v_{\tau,\beta},m_{\tau,\beta},\Lambda_{\tau,\beta})$ or, in case $\beta=0$ (stationary solutions), 
simply $(v_{\tau},m_{\tau},\Lambda_{\tau})$. 
We also notice that, if $(v,m,\Lambda)$ solves \eqref{eq:scalar}, then $u=\nu v$, $m$, 
$\lambda=\nu^2\Lambda$ solve
\[
\begin{cases}
 -\nu u''+\tfrac12|u'|^2+\xi u'+\lambda=\kappa R(1-\cos x),\quad u(0)=0,\\[2pt]
 -\nu m''-\bigl((u'+\xi)m\bigr)'=0,\quad m>0,\quad \ip m\,dx=1,
\end{cases}
 \qquad\text{ where } \tau=\frac{\kappa R}{\nu^2},\quad \beta=\frac\xi\nu .
\]

We have the following.

\begin{lemma}\label{lem:scalar}
For every $\tau\ge0$ and $\beta\in\R$, \eqref{eq:scalar} has a unique
solution, smooth, and $(\tau,\beta)\mapsto(v_{\tau,\beta},m_{\tau,\beta},\Lambda_{\tau,\beta})$
is smooth with values in $C^{k}(\T)^2\times\R$, for every $k$. Moreover
\[
 \Lambda_{\tau,-\beta}=\Lambda_{\tau,\beta},\qquad
 v_{\tau,-\beta}(x)=v_{\tau,\beta}(-x),\qquad
 m_{\tau,-\beta}(x)=m_{\tau,\beta}(x)=m_{\tau,\beta}(-x).
\]
For $\beta=0$, $v_\tau$ is even and $m_\tau=e^{-v_\tau}/\ip e^{-v_\tau}$.
At $\tau=0$, $v_{0,\beta}=0$, $m_{0,\beta}=(2\pi)^{-1}$, $\Lambda_{0,\beta}=0$.
\end{lemma}

\begin{proof}
The Hopf--Cole substitution $v=-2\log\psi$ transforms the first (uncoupled) equation of
\eqref{eq:scalar} into the eigenvalue problem
\[
 L_{\tau,\beta}\psi:=-2\psi''+2\beta\psi'+\tau(1-\cos x)\psi=\Lambda\psi\quad\text{ in }\T,
 \qquad \psi>0,\quad \psi(0)=1;
\]
by the Krein--Rutman theorem (see e.g. \cite[Section~I.7]{Hess}), this problem has a unique solution 
$\Lambda=\Lambda_{\tau,\beta}$, the unique principal eigenvalue of $L_{\tau,\beta}$, and 
$\psi=\psi_{\tau,\beta}$, its unique normalized positive eigenfunction, smooth. 
Since the potential $\tau(1-\cos x)$ is even, the reflection $x\mapsto-x$ transforms 
$L_{\tau,\beta}$ into $L_{\tau,-\beta}$, whence $\Lambda_{\tau,-\beta}=\Lambda_{\tau,\beta}$ and 
$\psi_{\tau,-\beta}(x)=\psi_{\tau,\beta}(-x)$. 
We deduce that the HJB equation has exactly one solution, with components 
$\Lambda = \Lambda_{\tau,\beta}$ and $v = v_{\tau,\beta} = -2\log\psi_{\tau,\beta}$, and that 
$v_{\tau,-\beta}(x)=v_{\tau,\beta}(-x)$ (in particular, for $\beta=0$, $v_{\tau}$ is even).

To solve the density equation, we notice that $L_{\tau,-\beta}$ is the formal adjoint of 
$L_{\tau,\beta}$, and we use its eigenfunction $\psi_{\tau,-\beta}$. Following 
\cite[Section~5.4]{ChetriteTouchette}, we set
\[
m=\frac{\psi_{\tau,\beta}\cdot\psi_{\tau,-\beta}}{\ip\psi_{\tau,\beta}\psi_{\tau,-\beta}\,dx}
\]
which is even, invariant under the change of $\beta$ into $-\beta$, and satisfies $\ip m=1$ 
(for $\beta=0$, it reduces to $e^{-v_\tau}/\ip e^{-v_\tau}$). A direct computation gives
\[
\bigl(m'+(v_{\tau,\beta}'+\beta)m\bigr)'=
\frac{\psi_{\tau,-\beta}\cdot L_{\tau,\beta}\psi_{\tau,\beta}-\psi_{\tau,\beta}\cdot L_{\tau,-\beta}
\psi_{\tau,-\beta}}{2\ip\psi_{\tau,\beta}\psi_{\tau,-\beta}\,dx}=0,
\]
thus $m$ solves the density equation, and the flux $J:=m'+(v_{\tau,\beta}'+\beta)m$ is constant. 
If $\rho=mf$ is another solution, then $(mf'+Jf)'=0$; testing with $f$ gives
$\ip m|f'|^2=0$, whence $f=1$. Then $m$ is the unique solution. In particular, 
the case $\tau=0$ is now explicit.

We are left to check smooth dependence on the parameters: since the principal eigenvalue is 
simple and isolated, this is a standard application of the implicit function theorem and Fredholm's alternative. 
\end{proof}

\begin{remark}\label{rem:flux-scalar}
Let $(v,m,\Lambda)$ solve \eqref{eq:scalar}. As we noticed
in the proof of Lemma~\ref{lem:scalar}, the flux of $m$ is constant; arguing
as for \eqref{eq:flux}, we obtain
\begin{equation}\label{eq:flux-scalar}
 m'+(v'+\beta)m=J,\qquad\text{where }\ J=J_{\tau,\beta}:=\frac{2\pi\beta}{\ip m^{-1}\,dx}
\end{equation}
(with a slight abuse of notation: after the above scaling, the flux in
\eqref{eq:flux} corresponds to $\nu J_{\tau,\beta}$). As a consequence, for
every $f\in C^2(\T)$,
\begin{equation}\label{eq:dual}
 m\bigl(-f''+(v'+\beta)f'\bigr)=-(mf')'+Jf',\qquad\text{whence }
 \ip m\bigl(-f''+(v'+\beta)f'\bigr)\,dx=0 .
\end{equation}
For $\beta=0$ we have $J=0$, thus $m'=-v' m$ and
$-(mf')'=m(-f''+v' f')$.
\end{remark}

Lemma \ref{lem:scalar} implies the following fixed point formulation of \eqref{eq:profiles}. 

\begin{lemma}\label{lem:reduction}
Let $(U,M,\bar\lambda)$ solve \eqref{eq:profiles}, let $R=R[M]$, and let $q\in\T$
satisfy $C[M]=R\cos q$, $S[M]=R\sin q$ ($q$ arbitrary if $R=0$). Then, with $\tau=\kappa R/\nu^2$ and $\beta=\xi/\nu$,
\begin{equation}\label{eq:reduction-xi}
\begin{gathered}
 U(\cdot+q)-U(q)=\nu v_{\tau,\beta},\qquad M(\cdot+q)=m_{\tau,\beta},\\
 \bar\lambda=\nu^2\Lambda_{\tau,\beta}-\kappa R,\qquad
 R=A\!\left(\frac{\kappa R}{\nu^2},\frac{\xi}{\nu}\right),
\end{gathered}
\end{equation}
where 
\begin{equation}\label{eq:def_respo_A}
A(\tau,\beta):=C[m_{\tau,\beta}],\qquad\text{ and }S[m_{\tau,\beta}]=0
\end{equation}
since $m_{\tau,\beta}$ is even; moreover, $A$ is even in $\beta$. Conversely, every $R\ge0$ with $R=A(\kappa R/\nu^2,\xi/\nu)$ produces in
this way a solution of \eqref{eq:profiles} with $C[M]=R$, $S[M]=0$. In particular,
$R=0$ gives exactly the incoherent solution.
\end{lemma}

\begin{proof}
Since $K[M](x)=-\kappa R\cos(x-q)$, the translated pair solves
\eqref{eq:scalar} after the above scaling, with $\Lambda=(\bar\lambda+\kappa R)/\nu^2$, and
Lemma~\ref{lem:scalar} identifies it. The self-consistency
$R=C[M(\cdot+q)]=C[m_{\tau,\beta}]$ is the last identity. The converse is the same
computation read backwards. Finally, $A$ is even in $\beta$ because $m_{\tau,-\beta}=m_{\tau,\beta}$.
\end{proof}

In the following we write $A(\tau)=A(\tau,0)$. After the above reduction, most results in 
Theorem \ref{thm:waves} will follow by a detailed analysis of the \emph{fixed-point map} 
$(\tau,\beta)\mapsto A(\tau,\beta)$, defined in \eqref{eq:def_respo_A}, and its derivatives 
with respect to $\tau$. In this direction, it is convenient to introduce the function $g_{\tau,\beta}$, 
defined through the following result.

\begin{lemma}\label{lem:corr-g}
Let $\tau\ge0$ and $\beta\in\R$, and write $m=m_{\tau,\beta}$, $A=A(\tau,\beta)$ and 
$J=J_{\tau,\beta}$ (recall that $J_{\tau,\beta}$ is defined in \eqref{eq:flux-scalar}). Then the problem
\begin{equation}\label{eq:corr-g-drift}
 -(mg')'+\frac{J^2}{m}\left(g-\frac{\ip g\,m^{-1}\,dx}{\ip m^{-1}\,dx}\right)=m\,(A-\cos x),
 \qquad \ip m g\,dx=0,
\end{equation}
has a unique solution $g=g_{\tau,\beta}\in H^1(\T)$, which is smooth and even. As usual, we write 
$g_\tau$ in case $\beta=0$, when \eqref{eq:corr-g-drift} reduces to
\begin{equation}\label{eq:corr-g}
 -(m_\tau g_\tau')'=m_\tau\,(A(\tau)-\cos x),\qquad \ip m_\tau g_\tau\,dx=0 .
\end{equation}
\end{lemma}

\begin{proof}
Since $m$ is bounded above and below, we can apply the Lax--Milgram theorem in the space of the 
functions $f\in H^1(\T)$ with $\ip mf\,dx=0$. Indeed, the bilinear form associated to the left 
hand side of \eqref{eq:corr-g-drift},
\[
 (g,f)\mapsto\ip mg'f'\,dx+J^2\left(\ip\frac{gf}{m}\,dx
 -\frac{\ip g\,m^{-1}\,dx\ip f\,m^{-1}\,dx}{\ip m^{-1}\,dx}\right),
\]
is symmetric, and it is coercive thanks to its first term and the Poincar\'e--Wirtinger inequality, the 
second one being nonnegative for $f=g$ by the Cauchy--Schwarz inequality. Furthermore, both sides 
of the equation vanish when tested with $f=1$ (recall that $\ip m\cos x\,dx=A$), thus the equation 
holds when tested with any $f\in H^1(\T)$. Finally, $g$ is smooth by elliptic regularity and, since 
$m$ is even, $g(-x)$ is again a solution; by uniqueness, $g$ is even.
\end{proof}

As a consequence, we can write the derivatives of $v_{\tau,\beta}$, $m_{\tau,\beta}$, $\Lambda_{\tau,\beta}$ and $A$ with respect to $\tau$ in terms of $g_{\tau,\beta}$.

\begin{lemma}\label{lem:resp}
The map $A$ is smooth on $[0,\infty)\times\R$ and even in $\beta$, $A(0,\beta)=0$, and $A<1$. 
Moreover, writing $v,m,\Lambda,g,J$ and $A$ for $v_{\tau,\beta},m_{\tau,\beta},\Lambda_{\tau,\beta},
g_{\tau,\beta},J_{\tau,\beta}$ and $A(\tau,\beta)$, and
\[
 \langle g\rangle:=\frac{\ip g\,m^{-1}\,dx}{\ip m^{-1}\,dx},
\]
we have
\begin{align}
 &\partial_\tau m=-g\,m,\qquad
  \partial_\tau v'=g'+\frac Jm\,\bigl(g-\langle g\rangle\bigr),\qquad
  \partial_\tau\Lambda=1-A,\label{eq:dv-drift}\\
 &\partial_\tau A=\ip m|\partial_\tau v'|^2\,dx
   =\ip m|g'|^2\,dx+J^2\ip\frac{(g-\langle g\rangle)^2}{m}\,dx>0,\label{eq:A1-drift}\\
 &\partial_\tau^2A=-3\ip m\,g\,|\partial_\tau v'|^2\,dx
   =-3\left[\ip m\,g\,|g'|^2\,dx+J^2\ip\frac{g\,(g-\langle g\rangle)^2}{m}\,dx\right],
   \label{eq:A2-drift}
\end{align}
and
\begin{equation}\label{eq:A1-zero}
 \partial_\tau A(0,\beta)=\frac1{2(1+\beta^2)} .
\end{equation}
In particular, for $\beta=0$,
\begin{align}
 \partial_\tau m_\tau&=-g_\tau m_\tau,\qquad
 \partial_\tau v_\tau=g_\tau-g_\tau(0),\label{eq:dv}\\
 A'(\tau)&=\ip m_\tau|g_\tau'|^2\,dx>0,\qquad A'(0)=\tfrac12,\label{eq:A1}\\
 A''(\tau)&=-3\ip m_\tau\,g_\tau\,|g_\tau'|^2\,dx .\label{eq:A2}
\end{align}
\end{lemma}

\begin{proof}
The first properties follow from Lemmas~\ref{lem:scalar} and \ref{lem:reduction}; in particular, 
$A<1$ since $m>0$. We write
\[
a=\partial_\tau v,\qquad n=\partial_\tau m,
\]
smooth by Lemma~\ref{lem:scalar}. Differentiating \eqref{eq:scalar} with
respect to $\tau$ we have
\[
\begin{cases}
 -a''+(v'+\beta)a'+\partial_\tau \Lambda=1-\cos x,\\[2pt]
 -n''-\bigl((v'+\beta)n\bigr)'=(ma')' ,\qquad \ip n\,dx=0 .
\end{cases}
\]
Multiplying the first equation by $m$ and integrating, \eqref{eq:dual} gives $\partial_\tau \Lambda=1-A$; hence
\begin{equation}\label{eq:lin-drift}
 -a''+(v'+\beta)a'=A-\cos x .
\end{equation}
On the other hand, testing the equation of $n$ with any $f\in C^2(\T)$,  we have
\begin{equation}\label{eq:dual-n}
 \ip n\bigl(-f''+(v'+\beta)f'\bigr)\,dx=-\ip m\,a'f'\,dx .
\end{equation}
Choosing $f=a$, and recalling \eqref{eq:lin-drift} and $\ip n=0$, we deduce that
\[
 \partial_\tau A=\ip n\cos x\,dx=-\ip n\,(A-\cos x)\,dx=\ip m|a'|^2\,dx ,
\]
which is positive, since $a'\equiv0$ in \eqref{eq:lin-drift} would force
$\cos x\equiv A$. Let now $b=\partial_\tau a$. Since $\partial_\tau v'=a'$,
differentiating \eqref{eq:lin-drift} we obtain
\begin{equation}\label{eq:lin2-drift}
 -b''+(v'+\beta)b'=\partial_\tau A-|a'|^2,
\end{equation}
and \eqref{eq:dual-n} with $f=b$ yields $\ip m\,a'b'=\ip n|a'|^2$. Then
\[
 \partial_\tau^2A=\ip n|a'|^2\,dx+2\ip m\,a'b'\,dx=3\ip n|a'|^2\,dx .
\]

We are left to write $n$ and $a'$ in terms of the function $g$. To this aim, let 
\[
\gamma=-n/m, \qquad\text{ so that }\ip m\gamma\,dx=0. 
\]
We are going to show that $\gamma = g$. Indeed, by \eqref{eq:flux-scalar},
$v'+\beta=(J-m')/m$. Differentiating with respect to $\tau$, and using $n=-\gamma m$,
\[
 a'=\gamma'+\frac{\partial_\tau J+J\gamma}{m}.
\]
Integrating on $\T$ gives $\partial_\tau J\ip m^{-1}+J\ip\gamma\,m^{-1}=0$,
that is $\partial_\tau J=-J\langle\gamma\rangle$, whence
\[
 a'=\gamma'+\frac Jm\,\bigl(\gamma-\langle\gamma\rangle\bigr).
\]
On the other hand, by \eqref{eq:dual} with $f=a$, equation \eqref{eq:lin-drift} can be written as 
$-(ma')'+Ja'=m(A-\cos x)$; inserting the above expression of $a'$, we see that $\gamma$ solves 
\eqref{eq:corr-g-drift}. By Lemma~\ref{lem:corr-g}, $\gamma=g$, and \eqref{eq:dv-drift} follows. Moreover, 
since $\langle g\rangle$ is constant, we have
$\ip g'(g-\langle g\rangle)=\ip g\,g'(g-\langle g\rangle)=0$; then, expanding the square in 
$|a'|^2=|g'+J(g-\langle g\rangle)/m|^2$, we obtain the last equalities in \eqref{eq:A1-drift} and 
\eqref{eq:A2-drift}.

At $\tau=0$ we have $m=(2\pi)^{-1}$, $A=0$ (Lemma~\ref{lem:scalar}) and $J=\beta/(2\pi)$; then 
\eqref{eq:corr-g-drift} reads $-g''+\beta^2g=-\cos x$, with $\ip g\,dx=0$, whence, using also \eqref{eq:dv-drift},
\begin{equation}\label{eq:g-zero}
 g_{0,\beta}=-\frac{\cos x}{1+\beta^2},\qquad
 \partial_\tau v_{\tau,\beta}'\big|_{\tau=0}=\frac{\sin x-\beta\cos x}{1+\beta^2} .
\end{equation}
Then \eqref{eq:A1-drift} gives
$\partial_\tau A(0,\beta)=\frac1{2\pi(1+\beta^2)^2}\ip|\sin x-\beta\cos x|^2\,dx=\frac1{2(1+\beta^2)}$.

Finally, let $\beta=0$. Then $J=0$, thus \eqref{eq:corr-g-drift} is \eqref{eq:corr-g} and 
$a'=g'$; since $a(0)=0$, we obtain \eqref{eq:dv}, \eqref{eq:A1} (recall \eqref{eq:A1-zero}) and \eqref{eq:A2}. We also notice 
that, in this case, $b$ and $\partial_\tau g_\tau$ differ by a constant, whence \eqref{eq:lin2-drift} and 
Remark~\ref{rem:flux-scalar} yield
\begin{equation}\label{eq:b-zero}
 -\bigl(m_\tau(\partial_\tau g_\tau)'\bigr)'=m_\tau\bigl(A'(\tau)-|g_\tau'|^2\bigr),
\end{equation}
which will be used below.
\end{proof}

\subsection{Concavity of the fixed-point map}\label{subsec:concavity}

In view of \eqref{eq:A2}, the fixed-point map of the stationary problem is strictly concave if 
and only if a suitable integral involving the function $g$ is positive. This is the content of the 
following result, which is the key step in the proof of Theorem~\ref{thm:waves}.

\begin{proposition}\label{prop:concavity}
$A''(\tau)<0$ for every $\tau>0$. Equivalently,
$\ip m_\tau\,g_\tau\,|g_\tau'|^2\,dx>0$ for $\tau>0$.
\end{proposition}

As we already mentioned, the proof of Proposition~\ref{prop:concavity} follows from an adaptation of the arguments in 
\cite{munoz2026uniqueness}. We report all the details in Appendix~\ref{sec:concavity}, for the sake of completeness. We also notice 
that, in the case $\beta=0$, the scalar problem \eqref{eq:scalar} can be written in terms of Mathieu 
functions; then Proposition~\ref{prop:concavity} can be equivalently stated as a property of the 
corresponding characteristic value.

\begin{remark}\label{rem:mathieu}
For $\tau\ge0$, the function $w(z)=\psi_\tau(2z)$, where $\psi_\tau=e^{-v_\tau/2}$ is the positive 
eigenfunction in the proof of Lemma~\ref{lem:scalar}, is positive, $\pi$-periodic and solves the 
Mathieu equation
\[
 w''+(a-2q\cos2z)\,w=0,\qquad\text{where }\ a=2(\Lambda_\tau-\tau),\quad q=-\tau .
\]
Since, among the $\pi$-periodic solutions of the Mathieu equation, only $\mathrm{ce}_0(\cdot,q)$ has no 
zeros, we have $a=a_0(-\tau)$, where $a_0$ denotes the lowest characteristic value 
(see \cite[Chapter~16]{Erdelyi3}, \cite{MeixnerSchafke}). Recalling that $a_0$ is even, we deduce that
\[
 \Lambda_\tau=\tau+\tfrac12\,a_0(\tau),
\]
whence, by \eqref{eq:dv-drift}, $A(\tau)=-\frac12\,a_0'(\tau)$. As a consequence, Lemma~\ref{lem:resp} 
and Proposition~\ref{prop:concavity} yield, for every $q>0$,
\begin{equation}\label{eq:mathieu}
 -2<a_0'(q)<0,\qquad a_0''(q)<0,\qquad a_0'''(q)>0
\end{equation}
(and $a_0'''(q)<0$ for $q<0$, since $a_0'''$ is odd). While the first two properties are elementary, 
we are not aware of a reference for the last one, which in particular says that $a_0'$ is strictly 
convex on $(0,\infty)$. Notice that the power series expansion of $a_0$ at $q=0$ only gives 
$a_0'''(q)=\frac{21}{16}\,q+O(q^3)$, in agreement with Lemma~\ref{lem:resp-zero} below, and that 
it converges only for $|q|$ smaller than about $1.47$ (see \cite{MeixnerSchafke}); on the other hand, 
the asymptotic expansion as $q\to\infty$ gives $a_0'''(q)\sim\frac34\,q^{-5/2}$ 
(see again \cite[Chapter~16]{Erdelyi3}).
\end{remark}

To deal with traveling waves, we extend the above concavity property to $\beta\ne0$ small, 
by a perturbation argument. We first need the following.

\begin{lemma}\label{lem:resp-zero}
For every $\beta\in\R$,
\[
 \partial_\tau^2A(0,\beta)=0,
\qquad\text{ and }\qquad
\partial_\tau^3A(0,0)=A'''(0)<0.
\]
\end{lemma}

\begin{proof}
By \eqref{eq:A2-drift} and \eqref{eq:g-zero} (recall that $m_{0,\beta}=(2\pi)^{-1}$),
\[
 \partial_\tau^2A(0,\beta)=\frac{3}{2\pi(1+\beta^2)^3}\ip\cos x\,(\sin x-\beta\cos x)^2\,dx=0,
\]
since $(\sin x-\beta\cos x)^2$ is a combination of $1$, $\cos2x$ and $\sin2x$.

Let now $\beta=0$, and write $m,g$ for $m_\tau,g_\tau$, and $b=\partial_\tau g_\tau$. 
Differentiating \eqref{eq:A2} and using \eqref{eq:dv},
\[
 A'''=-3\ip\Bigl(-m\,g^2|g'|^2+m\,b\,|g'|^2
   +2m\,g\,g'b'\Bigr)\,dx .
\]
At $\tau=0$ we have $m=(2\pi)^{-1}$, $g=-\cos x$ and $A'=\frac12$ (see \eqref{eq:g-zero} and 
\eqref{eq:A1}), thus \eqref{eq:b-zero} reads $-b''=\frac12-\sin^2x=\frac12\cos2x$. Moreover, 
differentiating the normalization $\ip m_\tau g_\tau=0$ we obtain $\ip m\,b=\ip m\,g^2$, that is 
$\ip b\,dx=\ip\cos^2x\,dx=\pi$ at $\tau=0$; hence $b=\frac18\cos2x+\frac12$. Since $2g\,g'b'=\frac14\sin^22x=\cos^2x\sin^2x=g^2|g'|^2$,
\[
 A'''(0)=-\frac3{2\pi}\ip b\,|g'|^2\,dx =-\frac3{2\pi}\ip\bigl(\tfrac18\cos2x+\tfrac12\bigr)
 \bigl(-\tfrac12\cos2x+\tfrac12\bigr)\,dx
   =-\frac{21}{32}.\qedhere
\]
\end{proof}

As a consequence, we obtain the strict concavity of $A(\cdot,\beta)$ on bounded intervals, 
provided $\beta$ is small enough.

\begin{proposition}\label{prop:concavity-drift}
For every $\tau_{\max}>0$ there is $\beta_0>0$ such that
\[
 \partial_\tau^2A(\tau,\beta)<0\qquad\text{for } 0<\tau\le\tau_{\max},\ |\beta|\le\beta_0 .
\]
\end{proposition}

\begin{proof}
By Lemma~\ref{lem:resp-zero}, $\partial_\tau^2A(0,\beta)=0$ for every
$\beta$, hence 
\[
 \partial_\tau^2A(\tau,\beta)=\tau\int_0^1\partial_\tau^3A(s\tau,\beta)\,ds =:\tau\, h(\tau,\beta),
\]
where $h$ is a continuous function of $(\tau,\beta)\in[0,\infty)\times\R$. 
Since
\[
h(\tau,0) = 
\begin{cases}
A''(\tau)/\tau & \tau>0,\\
A'''(0) & \tau=0,
\end{cases}
\] 
Proposition~\ref{prop:concavity} and  Lemma~\ref{lem:resp-zero} imply that $h<0$ on the compact set 
$[0,\tau_{\max}]\times\{0\}$ and, by continuity, on 
$[0,\tau_{\max}]\times[-\beta_0,\beta_0]$ for some $\beta_0>0$.
\end{proof}

We are now in a position to prove Theorem~\ref{thm:waves}. Since the first part of the 
argument will be used also in the proof of Theorem~\ref{thm:classification}, we do not 
assume $\kappa>\kappa_c$ from the beginning.

\begin{proof}[Proof of Theorem~\ref{thm:waves}]
Let $\kappa>0$, let $\beta_0>0$ be given by Proposition~\ref{prop:concavity-drift} with 
$\tau_{\max}=\kappa/\nu^2$, and let $|\xi|\le\nu\beta_0$. By Lemma~\ref{lem:reduction}, 
writing $\beta=\xi/\nu$, we are left to study the positive zeros of
\begin{equation}\label{eq:Gamma}
 \Gamma(R):=R-A\!\left(\frac{\kappa R}{\nu^2},\beta\right),\qquad R\ge0 .
\end{equation}
Since $A<1$, they lie in $(0,1)$, and $\Gamma(1)>0$. Moreover, by 
Lemma~\ref{lem:resp} (recall \eqref{eq:A1-zero} and \eqref{eq:kappac}),
\[
 \Gamma(0)=0,\qquad \Gamma'(0)=1-\frac{\kappa}{\kappa_c(\xi)} ,
\]
while, if $0\le R\le1$, then $\kappa R/\nu^2\le\tau_{\max}$; thus $\Gamma$ is strictly 
convex on $[0,1]$ by Proposition~\ref{prop:concavity-drift}, and in particular
\begin{equation}\label{eq:tangent}
 \Gamma(R)>\Gamma'(0)R,\qquad 0<R\le1 .
\end{equation}
If $\kappa\le\kappa_c(\xi)$, then \eqref{eq:tangent} gives $\Gamma>0$ on $(0,1]$: there 
are no positive zeros, and the incoherent solution is the only solution of 
\eqref{eq:profiles}. If $\kappa>\kappa_c(\xi)$, then $\Gamma<0$ near $0^+$ and 
$\Gamma(1)>0$, thus $\Gamma$ has a zero $R_\xi\in(0,1)$; by strict convexity, it is the 
only positive one, and
\begin{equation}\label{eq:radial-gap-xi}
 \delta_\xi:=1-\frac{\kappa}{\nu^2}\,
  \partial_\tau A\!\left(\frac{\kappa R_\xi}{\nu^2},\frac\xi\nu\right)
  =\Gamma'(R_\xi)>\frac{\Gamma(R_\xi)-\Gamma(0)}{R_\xi}=0 .
\end{equation}
In such a case, we define $(U_\xi,M_\xi,\bar\lambda_\xi)$ through
\eqref{eq:reduction-xi}, with $R=R_\xi$ and $q=0$: by Lemma~\ref{lem:reduction}, it is a 
synchronized solution of \eqref{eq:profiles}, with $C[M_\xi]=R_\xi$ and $S[M_\xi]=0$, 
and any other synchronized solution is a translate of it.

Let now $\kappa>\kappa_c$, and let $\xi$ satisfy \eqref{eq:speeds}. Then (i) has just been 
proved. As for (ii), the function $(R,\xi)\mapsto R-A(\kappa R/\nu^2,\xi/\nu)$ is smooth
and its derivative with respect to $R$ at $(R_\xi,\xi)$ is $\delta_\xi>0$; by the implicit
function theorem and the uniqueness of $R_\xi$, the map $\xi\mapsto R_\xi$ is smooth, and
Lemma~\ref{lem:scalar} gives the smoothness of the profiles. Since $A$ is even
in $\beta$, $R_{-\xi}=R_\xi$, and the symmetries in (iii) are those of
Lemma~\ref{lem:scalar}. Finally, for $\xi=0$ we have $R_0=C[M_0]$; let $\tau_0=\kappa R_0/\nu^2$. 
Then $M_0=m_{\tau_0}$ and $A(\tau_0)=R_0$, whence $g_0=g_{\tau_0}$ by \eqref{eq:corr-g}, 
and \eqref{eq:A1} gives $\delta_0=1-\frac{\kappa}{\nu^2}A'(\tau_0)=1-\frac{\kappa}{\nu^2}\ip M_0|g_0'|^2\,dx$; then (iv) 
follows from \eqref{eq:radial-gap-xi}.
\end{proof}

We also give the proof of Theorem~\ref{thm:classification}; recall that, in the first part of the previous proof, $\kappa>0$ is arbitrary.

\begin{proof}[Proof of Theorem~\ref{thm:classification}]
Part (b) has been obtained in the first part of the proof of Theorem~\ref{thm:waves}, 
which holds for every $\kappa>0$, with $\beta_0$ given by 
Proposition~\ref{prop:concavity-drift} with $\tau_{\max}=\kappa/\nu^2$. Taking $\xi=0$ 
there, we also obtain the alternative in part (a); moreover, in case 
$\kappa>\kappa_c$, the nondegeneracy of $(U_0,M_0,\bar\lambda_0)$ is 
Theorem~\ref{thm:waves}(iv). We are left to check the properties of the map 
$\kappa\mapsto R_0(\kappa)$, $\kappa>\kappa_c$. To this aim, recall that 
$R_0(\kappa)$ is the unique positive zero of
\[
 \Gamma_\kappa(R)=R-A\!\left(\frac{\kappa R}{\nu^2}\right),\qquad R\ge0 ,
\]
that $\Gamma_\kappa(0)=0$, $\Gamma_\kappa'(0)=1-\kappa/\kappa_c$, and that 
$\Gamma_\kappa$ is strictly convex on $[0,\infty)$, for every $\kappa>0$, by 
Proposition~\ref{prop:concavity}. The map $(\kappa,R)\mapsto\Gamma_\kappa(R)$ is smooth and
$\partial_R\Gamma_\kappa(R_0(\kappa))=\delta_0>0$; then the implicit function theorem and
the uniqueness of $R_0(\kappa)$ give the smoothness of $\kappa\mapsto R_0(\kappa)$ on 
$(\kappa_c,\infty)$. Moreover, differentiating the identity $\Gamma_\kappa(R_0(\kappa))=0$ with respect to $\kappa$, we obtain $\delta_0\,R_0'(\kappa)=\nu^{-2}R_0(\kappa)\,A'(\tau_0)>0$, where $\tau_0=\kappa R_0(\kappa)/\nu^2$, whence $R_0$ is strictly increasing. Finally, fix $\varepsilon>0$. Since $\Gamma_{\kappa_c}'(0)=0$, strict 
convexity gives $\Gamma_{\kappa_c}(\varepsilon)>0$,
hence $\Gamma_\kappa(\varepsilon)>0$ for $\kappa$ close to $\kappa_c$. On the other hand,
$\Gamma_\kappa<0$ on $(0,R_0(\kappa))$, again by strict convexity; we deduce that 
$R_0(\kappa)<\varepsilon$ for such $\kappa$.
\end{proof}

\begin{remark}\label{rem:kappac}
By \eqref{eq:A1-zero}, the function $\Gamma$ in \eqref{eq:Gamma} satisfies
$\Gamma'(0)=1-\kappa/\kappa_c(\xi)$ for every $\xi\in\R$: traveling waves of
speed $\xi$ bifurcate from the incoherent solution at
$\kappa=\kappa_c(\xi)$. We want to stress here that the restriction to small speeds in
Proposition~\ref{prop:concavity-drift} cannot be removed altogether.
Expanding in $\tau$ the Fourier coefficients of the eigenfunction
$\psi=e^{-v_{\tau,\beta}/2}$ (we omit the computation) one finds
\[
 \partial_\tau^3A(0,\beta)=\frac{3\,(5\beta^2-7)}{8(1+\beta^2)^3(4+\beta^2)} .
\]
Hence, if $\xi^2>7\nu^2/5$, $A(\cdot,\xi/\nu)$ is strictly convex near
$\tau=0$ and the bifurcation is subcritical: for $\kappa<\kappa_c(\xi)$ close
to $\kappa_c(\xi)$, $\Gamma$ is positive near $0^+$ and at $R=1$, and negative
somewhere in between. Then \eqref{eq:profiles} has at least two solutions
with $R[M]>0$ which are not translates of each other, although
$\kappa<\kappa_c(\xi)$.
\end{remark}

\subsection{Momentum and coercivity}\label{subsec:coercivity}

In the last part of this section, we address some dynamical properties of traveling waves, which 
will be used in the following. Accordingly, in the rest of the section $\kappa>\kappa_c$ is fixed, 
as well as $\xi$ in a suitable open interval containing $0$, so that the waves 
$(U_\xi,M_\xi,\bar\lambda_\xi)$ are defined and depend smoothly on $\xi$, by
Theorem~\ref{thm:waves}.

To start with, we compute the momentum of traveling waves; indeed, recall that
the \emph{momentum} $P=\ip mu_x\,dx$ of a solution of \eqref{eq:mfg} is conserved, see e.g. Proposition~\ref{prop:momentum} ahead.

\begin{lemma}\label{lem:flux}
Let $(U_\xi,M_\xi,\bar\lambda_\xi)$ solve \eqref{eq:profiles}, let $J_\xi$ be its 
flux (see \eqref{eq:flux}). Then 
\begin{equation}\label{eq:momentum}
 P_\xi:=\ip M_\xi U_\xi'\,dx=2\pi J_\xi-\xi=\xi\left(\frac{4\pi^2}{\ip M_\xi^{-1}\,dx}-1\right).
\end{equation}
Moreover $\ip M_\xi^{-1}\,dx\ge4\pi^2$, with equality if and only if $M_\xi$ is
constant; in particular $|J_\xi|\le|\xi|/(2\pi)$.
\end{lemma}

\begin{proof}
Integrating the first identity in \eqref{eq:flux} gives
$2\pi J_\xi=P_\xi+\xi$, and \eqref{eq:momentum} follows from the second one. Finally,
$4\pi^2=\bigl(\ip M_\xi^{1/2}M_\xi^{-1/2}\bigr)^2\le\ip M_\xi\ip M_\xi^{-1}$ by
the Cauchy--Schwarz inequality, with equality if and only if $M_\xi$ is
constant.
\end{proof}

Next, we turn to some coercivity properties, related to the dynamical stability of traveling waves. 
We start from the stationary case. 
Let
\[
 \|f\|_\xi^2:=\ip M_\xi f^2\,dx,\qquad
 H^1_\xi:=\Bigl\{f\in H^1(\T):\ \ip M_\xi f\,dx=0\Bigr\},
\]
and, for $\theta,h\in H^1_\xi$,
\begin{equation}\label{eq:Qdef}
 \Q_\xi(\theta,h):=\|\theta'+\nu h'\|_\xi^2+\nu^2\|h'\|_\xi^2
   -\kappa\bigl(C[M_\xi h]^2+S[M_\xi h]^2\bigr).
\end{equation}
For $\xi=0$ we write $\|\cdot\|_0$ and $H^1_0$ (no boundary condition is
involved). In the next section, $\theta$ and $\eta=M_\xi h$ will be the
perturbations of the value function and of the density around a traveling
wave; note that $h\in H^1_\xi$ means $\ip\eta\,dx=0$. The form $\Q_0$
degenerates exactly in the direction of translations.

\begin{lemma}\label{lem:stationary-gap}
For every $\theta,h\in H^1_0$ we have $\Q_0(\theta,h)\ge0$, with equality if
and only if $(\theta,h)=c\,(U_0',M_0'/M_0)$ for some $c\in\R$.
\end{lemma}

\begin{proof}
Write $M=M_0$, $g=g_0$ and $\tau_0=\kappa R_0/\nu^2$ (recall that $g_0$ and $R_0$ are 
defined in Theorem~\ref{thm:waves}). By \eqref{eq:flux} with $\xi=0$, 
$U_0'=-\nu M'/M$, and $M'/M\in H^1_0$.
Differentiating in $x$ the first equation in \eqref{eq:profiles} (with $\xi=0$), where
$K[M]=-\kappa R_0\cos x$, we have $-\nu U_0'''+U_0'U_0''=\kappa R_0\sin x$,
that is $\bigl(M(M'/M)'\bigr)'=\tau_0M\sin x$. This and the equation of $g$
in Theorem~\ref{thm:waves} give, for every $f\in H^1_0$,
\begin{equation}\label{eq:weak-gr}
 \ip Mg'f'\,dx=-C[Mf],\qquad \ip M\Bigl(\frac{M'}M\Bigr)'f'\,dx=-\tau_0S[Mf].
\end{equation}
Taking $f=g$ and $f=M'/M$ in the second identity, and recalling
\eqref{eq:radial-gap},
\begin{equation}\label{eq:gr-norms}
 \ip Mg'\Bigl(\frac{M'}M\Bigr)'\,dx=0,\qquad
 \Bigl\|\Bigl(\frac{M'}M\Bigr)'\Bigr\|_0^2=\tau_0R_0,\qquad
 \|g'\|_0^2=\frac{\nu^2}{\kappa}\,(1-\delta_0) ;
\end{equation}
indeed $S[Mg]=0$ because $Mg$ is even, and, integrating by parts, $S[M']=-C[M]=-R_0$.

Let $h\in H^1_0$ and write
\[
 h=a\,\frac{M'}M+bg+k,\qquad\text{where }\ 
 a=\frac{\ip Mh'(M'/M)'\,dx}{\|(M'/M)'\|_0^2},\quad
 b=\frac{\ip Mh'g'\,dx}{\|g'\|_0^2}.
\]
Then $k\in H^1_0$, $\ip Mk'(M'/M)'=\ip Mk'g'=0$, and
$\|h'\|_0^2=a^2\|(M'/M)'\|_0^2+b^2\|g'\|_0^2+\|k'\|_0^2$. By
\eqref{eq:weak-gr} and \eqref{eq:gr-norms}, $C[Mh]=-b\|g'\|_0^2$ and
$S[Mh]=-aR_0 $. Since $\nu^2\tau_0=\kappa R_0$,
\[
 \nu^2\|h'\|_0^2-\kappa\bigl(C[Mh]^2+S[Mh]^2\bigr)
 =b^2\|g'\|_0^2\bigl(\nu^2-\kappa\|g'\|_0^2\bigr)+\nu^2\|k'\|_0^2
 =\nu^2\delta_0\,b^2\|g'\|_0^2+\nu^2\|k'\|_0^2 .
\]
Therefore
\[
 \Q_0(\theta,h)=\|\theta'+\nu h'\|_0^2+\nu^2\delta_0\,b^2\|g'\|_0^2
   +\nu^2\|k'\|_0^2\ \ge0 .
\]
Since $\delta_0>0$, $g'\not\equiv0$, and a function of $H^1_0$ with zero
derivative vanishes, equality holds if and only if $b=0$, $k=0$ and
$\theta=-\nu h$, that is $(\theta,h)=a\,(U_0',M_0'/M_0)$.
\end{proof}

The analogous coercivity property for $\xi\neq0$ is contained in the next theorem, which also collects all the properties of the traveling waves used in the rest of the paper.

\begin{theorem}\label{thm:travel-gap}
Let $\kappa>\kappa_c$. There are $\xi_{\max}>0$ and positive constants $C_0$, $\underline R$, $c_P$, $c_0$ such that the following holds, with
$I=[-\xi_{\max},\xi_{\max}]$.
\begin{enumerate}
\item[(a)] \emph{Profiles.} Every $\xi\in I$ satisfies \eqref{eq:speeds}; thus 
 Theorem~\ref{thm:waves} provides the synchronized solution 
 $(U_\xi,M_\xi,\bar\lambda_\xi)$ of \eqref{eq:profiles}, which is unique up to 
 translations and depends smoothly on $\xi\in I$. Moreover, for every
 $\xi,\xi'\in I$ and $s\in\R$
 \begin{gather}
  0<C_0^{-1}\le M_\xi\le C_0,\qquad
  \|U_\xi\|_{C^2}+\|M_\xi\|_{C^1}\le C_0,\qquad
  1>R[M_\xi]\ge \underline R>0,\label{eq:profile-bounds}\\[2pt]
  \ip\frac{|M_\xi'|^2}{M_\xi}\,dx\ge\frac{\underline R^2}{1-\underline R^2}>0,\label{eq:profile-bounds2}\\[2pt]
  \|M_\xi(\cdot+s)-M_{\xi'}\|_{L^2}+\|U_\xi'(\cdot+s)-U_{\xi'}'\|_{L^2}
   \le C_0\,(|s|+|\xi-\xi'|).\label{eq:noescape-translate-profiles}
 \end{gather}
\item[(b)] \emph{Flux and momentum.} The flux $J_\xi$ and the momentum $P_\xi$ of 
 $(U_\xi,M_\xi,\bar\lambda_\xi)$ satisfy \eqref{eq:flux} and \eqref{eq:momentum}, 
 $|J_\xi|\le|\xi|/(2\pi)$, and $\xi\mapsto P_\xi$ is
 smooth and odd on $I$, with
 \begin{equation}\label{eq:working}
  -c_P^{-1}\le P_\xi'\le-c_P ,\qquad \xi\in I .
 \end{equation}
\item[(c)] \emph{Coercivity.} For every $\xi\in I$ and every
 $\theta,h\in H^1_\xi$ such that
 \begin{equation}\label{eq:phase-constraint}
  \ip M_\xi\bigl(U_\xi'h-\theta'\bigr)\,dx=0,
 \end{equation}
 we have
 \begin{equation}\label{eq:travel-coercivity}
  \Q_\xi(\theta,h)\ge c_0\bigl(\|\theta'\|_\xi^2+\|h'\|_\xi^2\bigr).
 \end{equation}
\end{enumerate}
\end{theorem}

\begin{proof}
\emph{(a)} Let $\beta_0$ be as in Theorem~\ref{thm:waves} and fix
$0<\xi_{\max}<\nu\beta_0$ with $\kappa>\kappa_c(\xi_{\max})$. Then \eqref{eq:speeds} holds 
in a neighborhood of $I$, where Theorem~\ref{thm:waves}
gives the waves, their uniqueness, and their smoothness. As a consequence, 
\eqref{eq:profile-bounds} holds by continuity and compactness,
with $\underline R=\min_IR[M_\xi]$, and \eqref{eq:noescape-translate-profiles} holds with
\[
C_0 \ge \max_{\xi\in I}\bigl(\|M_\xi'\|_{L^2}+\|U_\xi''\|_{L^2}
   +\|\partial_\xi M_\xi\|_{L^2}+\|\partial_\xi U_\xi'\|_{L^2}\bigr).
\]
For \eqref{eq:profile-bounds2}, let $R=C[M_\xi]$. Integrating by parts,
$R=-\ip\sin x\,M_\xi'\,dx$, hence
\[
 R^2\le\ip\sin^2x\,M_\xi\,dx\ip\frac{|M_\xi'|^2}{M_\xi}\,dx
    \le(1-R^2)\ip\frac{|M_\xi'|^2}{M_\xi}\,dx ,
\]
because $\ip\cos^2x\,M_\xi\ge\bigl(\ip\cos x\,M_\xi\bigr)^2$; recall that
$R\ge \underline R$.

\emph{(b)} By Lemma~\ref{lem:flux} and (a), $P_\xi$ is smooth; it is odd
because $M_{-\xi}=M_\xi$. Moreover
\[
 P_0'=\frac{4\pi^2}{\ip M_0^{-1}\,dx}-1<0,
\]
since $M_0$ is not constant. By continuity, \eqref{eq:working} holds for some
$c_P>0$ after reducing $\xi_{\max}$.

\emph{(c)} We argue by contradiction. If (c) fails for every choice of
$\xi_{\max}$ and $c_0$, there are $\xi_n\to0$ and
$\theta_n,h_n\in H^1_{\xi_n}$ satisfying \eqref{eq:phase-constraint} with
$\xi=\xi_n$, such that
\[
 \|\theta_n'\|_{\xi_n}^2+\|h_n'\|_{\xi_n}^2=1,\qquad
 \Q_{\xi_n}(\theta_n,h_n)<\frac1n .
\]
A function $f\in H^1_{\xi_n}$ vanishes somewhere, thus
\[
 \|f\|_\infty\le\sqrt{2\pi}\,\|f'\|_{L^2}\le(2C_0\pi)^{1/2}\|f'\|_{\xi_n}
\]
by \eqref{eq:profile-bounds}. Hence $(\theta_n)_n$ and $(h_n)_n$ are bounded in
$H^1(\T)$ and, up to a subsequence, they converge weakly in $H^1(\T)$ and
uniformly to some $\theta$ and $h$. Since $M_{\xi_n}\to M_0$ and
$U_{\xi_n}'\to U_0'$ uniformly, we get $\theta,h\in H^1_0$,
\[
 \ip M_0\bigl(U_0'h-\theta'\bigr)\,dx=0,\qquad
 C[M_{\xi_n}h_n]\to C[M_0h],\qquad S[M_{\xi_n}h_n]\to S[M_0h].
\]
Moreover, if $f_n\rightharpoonup f$ weakly in $L^2(\T)$, then
$\|f\|_0\le\liminf_n\|f_n\|_{\xi_n}$, because
$\bigl|\|f_n\|_{\xi_n}^2-\|f_n\|_0^2\bigr|\le\|M_{\xi_n}-M_0\|_\infty\|f_n\|_{L^2}^2\to0$.
Therefore
\[
 \Q_0(\theta,h)\le\liminf_n\Q_{\xi_n}(\theta_n,h_n)\le0 .
\]
By Lemma~\ref{lem:stationary-gap}, $(\theta,h)=c\,(U_0',M_0'/M_0)$, and the
constraint gives
\[
 0=c\ip\bigl(U_0'M_0'-M_0U_0''\bigr)\,dx=2c\ip U_0'M_0'\,dx
  =-2c\,\nu\ip\frac{|M_0'|^2}{M_0}\,dx ,
\]
whence $c=0$ by \eqref{eq:profile-bounds2}. Thus $\theta=h=0$, the integral coefficients
$C[M_{\xi_n}h_n]$, $S[M_{\xi_n}h_n]$ tend to zero, and
\[
 \|\theta_n'+\nu h_n'\|_{\xi_n}^2+\nu^2\|h_n'\|_{\xi_n}^2
 <\frac1n+\kappa\bigl(C[M_{\xi_n}h_n]^2+S[M_{\xi_n}h_n]^2\bigr)\to0 .
\]
Then $\|h_n'\|_{\xi_n}\to0$ and
$\|\theta_n'\|_{\xi_n}\le\|\theta_n'+\nu h_n'\|_{\xi_n}+\nu\|h_n'\|_{\xi_n}\to0$,
which contradicts the normalization. Then, for a possibly smaller $\xi_{\max}$, we obtain the full theorem.
\end{proof}

From now on $\kappa>\kappa_c$ is fixed, together with $\xi_{\max}$, $I$ and the
constants of Theorem~\ref{thm:travel-gap}; all the constants below may depend
on them.

\section{Analysis of the time-dependent problem}\label{sec:tdp}

In this section $(u,m)$ is a given solution of
\eqref{eq:mfg}, with
\[
 (u,m)\in C^{1+\alpha/2,\,2+\alpha}([0,T]\times\T)^2,\qquad
 m>0,\qquad \ip m(t,x)\,dx=1,\qquad 0<\alpha<1.
\]
We assume the existence of $\xi, a, r$ such that $0 < r \le 1$,
\begin{equation}\label{eq:apriori}
\sup_{t \in [0,T]} \left(\|m(t, \cdot)-M_{\xi}(\cdot-\xi t - a)\|_{L^2}+\|u_x(t, \cdot)-U'_{\xi}(\cdot-\xi t - a)\|_{L^2}\right) \le r,
\end{equation}
that is, $(u,m)$ is close to a wave with speed $\xi$. Proposition \ref{prop:momentum} below shows that the momentum $P_* = \ip m u_x\,dx$ is time-independent. The main goal of this section is to compare $(u,m)$ with the family of traveling waves $(U_{\xi_*}(x-\xi_* t-a), M_{\xi_*}(x-\xi_* t-a))$, where $\xi_*$ is the uniquely determined wave speed matching $P_*=P_{\xi_*}\,  ( = \int M_{\xi_*}U'_{\xi_*})$. 

\begin{proposition}[Momentum conservation]\label{prop:momentum}
For a given solution $(u,m)$ of \eqref{eq:mfg},
\begin{equation}\label{eq:conserved}
 P_*=\ip m(t,x)u_x(t,x)\,dx
\end{equation}
is independent of $t$. 
\end{proposition}

\begin{proof}
Let $P_*(t) = \ip mu_x dx$. Then integration by parts gives
\[
\frac{d}{dt} P_*(t) =\ip(m_tu_x-m_xu_t)\,dx.
\]
Substituting the equations,
\begin{align*}
\frac{d}{dt} P_*(t)
 &=\ip\left[
    \bigl(\nu m_{xx}+\partial_x(m u_x)\bigr)u_x
    -m_x\left(-\nu u_{xx}+\frac12u_x^2-K[m]\right)\right]dx\\
 &=\ip\partial_x\left(\nu m_xu_x+\tfrac12m u_x^2\right)dx
    +\ip m_x K[m]\,dx
 =-\ip m \partial_x K[m]\,dx.
\end{align*}
Therefore,
\[
\frac{d}{dt} P_*(t)=-\kappa\ip m(C[m]\sin x-S[m]\cos x)\,dx
          =-\kappa(C[m]S[m]-S[m]C[m])=0.
\]
\end{proof}

\begin{lemma}\label{lem:xizeta} There exists $r_0 > 0$, $0<\xi_0\le\xi_{\max}$, such that, if $r \le r_0$ and $|\xi| \le \xi_0$, then there exists a unique $\xi_* \in I$ such that
\begin{equation}\label{eq:matched-speed}
P_{\xi_*}=P_*.
\end{equation}
Moreover,
\[
|\xi_*- \xi| \le C r
\]
for some $C$ depending only on $\xi_{\max}$.
\end{lemma}

\begin{proof} Since
\begin{align}\label{eq:123}
P_* & - P_{\xi} = \ip m(t,x+\xi t+a)u_x(t,x+\xi t+a)\,dx - \ip M_{\xi}(x)U'_{\xi}(x)\,dx \notag\\
& = \ip (m-M_{\xi})(u_x - U'_{\xi}) + (m-M_{\xi})U'_{\xi} + (u_x - U'_{\xi}) M_{\xi}\, dx,
\end{align}
by \eqref{eq:apriori} and \eqref{eq:profile-bounds} there exists a constant $C_1$ depending on $\xi_{\max}$ such that
\[
|P_* - P_{\xi}| \le C_1 r.
\]
Since $-c^{-1}_P \le P' \le -c_P$ by \eqref{eq:working} and $P_0 = 0$, we have $|P_{\xi}| \le c_P^{-1} |\xi|$. By the triangle inequality,
\begin{equation}
|P_*| \le C_1 r +  c_P^{-1} |\xi|.
\end{equation}
If $r \le |P_{\xi_{\max}}|/(2C_1)$ and $|\xi| \le c_P |P_{\xi_{\max}}|/2$, we get $|P_*| \le |P_{\xi_{\max}}|$ (notice that $P_{\xi_{\max}} = -P_{-\xi_{\max}}$), hence by monotonicity of 
$P$ on $I$ we have the existence of a unique $|\xi_*| \le \xi_{\max}$ such that \eqref{eq:matched-speed} holds.

The last inequality follows by observing that $|\xi_*-\xi|\le c_P^{-1}|P_*-P_{\xi}|$.
\end{proof}

\begin{lemma}\label{lem:phase} There exists $r_0, \xi_0 > 0$ such that, if $r \le r_0$ and $|\xi| \le \xi_0$, then there exists $\tilde q \in C^1([0,T])$ and $C > 0$ such that
\begin{equation}\label{eq:crucial-identity}
\ip\bigl[U_{\xi_*}'(x)(m(t,x + \xi t+a+\tilde q(t))-M_{\xi_*}(x))
              -M_{\xi_*}(x)(u_x(t,x + \xi t+a+\tilde  q(t))-U_{\xi_*}'(x))\bigr]\,dx = 0.
\end{equation}
and 
\begin{align}\label{eq:phase-residual}
 |\tilde  q(t)| \le C\left|\ip\bigl[U_{\xi_*}'(x)(m(t, x+\xi t+a)-M_{\xi_*}(x))
              -M_{\xi_*}(x)(u_x(t, x+\xi t+a)-U_{\xi_*}'(x))\bigr]\,dx \right|, \\
 \norm{u_x(t,\cdot + \xi t+a+ \tilde q(t))-U_{\xi_*}'}_{L^2}+\norm{m(t,\cdot + \xi t+a+ \tilde q(t))-M_{\xi_*}}_{L^2}\label{eq:boh}
       \le C r.
\end{align}
\end{lemma}

\begin{proof} Let $r \le r_0$ and $|\xi| \le \xi_0$ as in the previous lemma. We will possibly shrink $r_0$ below.
Fix $t \in [0,T]$ and let, for $p, \mu \in L^2(\T)$ and $q \in \T$,
\begin{align*}
F(p,\mu,q)
 &=\ip\bigl[U_{\xi_*}'(x)(\mu(x+q)-M_{\xi_*}(x))
              -M_{\xi_*}(x)(p(x+q)-U_{\xi_*}'(x))\bigr]\,dx \\
 &=\ip\bigl[U_{\xi_*}'(y-q)\mu(y)-M_{\xi_*}(y-q)p(y)\bigr]\,dy.
 \end{align*}
We want to show that there exists $\tilde q$ such that, for
\begin{align*}
p(x) = u_x(t, x+\xi t+a), \quad \mu(x) = m(t, x+\xi t+a),
\end{align*}
we have $F(p, \mu, \tilde  q) = 0$. Note that $F$ is continuously Fr\'echet differentiable, and $F(U_{\xi_*}',M_{\xi_*},0)=0$. Moreover,
\begin{equation}\label{eq:Fq}
 \partial_qF(p,\mu,q)=
 \ip\bigl[M_{\xi_*}'(y-q)p(y)-U_{\xi_*}''(y-q)\mu(y)\bigr]\,dy.
\end{equation}
By \eqref{eq:flux} we have $U_{\xi_*}'=J_{\xi_*}/M_{\xi_*}-\nu M_{\xi_*}'/M_{\xi_*}-{\xi_*}$.
Since $\ip M_{\xi_*}'\,dx=\ip M_{\xi_*}'/M_{\xi_*}\,dx=0$, we get, using \eqref{eq:profile-bounds2}
\begin{equation}\label{eq:1000}
\partial_qF(U_{\xi_*}',M_{\xi_*},0) =  2\ip U_{\xi_*}'M_{\xi_*}'\,dx=-2\nu\ip\frac{|M_{\xi_*}'|^2}{M_{\xi_*}}\,dx \le -\frac{2\nu \underline R^2}{1-\underline R^2} < 0.
\end{equation}
Hence, there exists a neighborhood $\mathcal U$ of $(U_{\xi_*}',M_{\xi_*})$ in $L^{2}(\T) \times L^{2}(\T)$ and a neighborhood $\mathcal V$ of $0$ such that, if $(p,\mu)$ belongs to $\mathcal U$, then $F(p, \mu, q) = 0$ as desired for a uniquely determined $\tilde q \in \mathcal V$. Note that 
\[
\|p-U'_{\xi_*}\| + \|\mu-M_{\xi_*}\| \le \|p-U'_{\xi}\| + \|\mu-M_{\xi}\| + \|M_{\xi_*}-M_{\xi}\| + \|U'_{\xi_*}-U'_{\xi}\| \le r + Cr
\]
by \eqref{eq:noescape-translate-profiles} and Lemma \ref{lem:xizeta}. Hence $(p,\mu)$ belongs to $\mathcal U$ for $r$ small enough. Moreover, the mean value theorem gives $|\tilde q| = |F(p,\mu,0)| |F_q(p,\mu,\delta)|^{-1}$, hence \eqref{eq:phase-residual} follows by \eqref{eq:1000} and the fact that $p,\mu$ is close to $(U_{\xi_*}',M_{\xi_*})$. It is important to stress that $\mathcal U$, $\mathcal V$ are uniform with respect to $\xi_*$. This is because \eqref{eq:1000} and \eqref{eq:profile-bounds}, and the fact that $F$ is affine in $(p,\mu)$.

The last inequality is a consequence of
\begin{align*}
\norm{u_x(\cdot + \xi t+a+\tilde q)-U_{\xi_*}'}_{L^2} & = \norm{u_x(\cdot + \xi t+a)-U_{\xi_*}'(\cdot - \tilde q)}_{L^2} \\
& \le \norm{u_x(\cdot + \xi t+a)-U_{\xi_*}'}_{L^2} + \norm{U_{\xi_*}'(\cdot - \tilde q)-U_{\xi_*}'}_{L^2} \\
& \le \norm{p-U_{\xi_*}'}_{L^2} + C_0|\tilde q|,
\end{align*}
and a similar estimate for $\norm{m(\cdot + \xi t+a+\tilde q)-M_{\xi_*}}_{L^2}$.

Note finally that $\tilde q$, now as a function of time, is of class $C^1$. Indeed, the map 
\[
t\mapsto F(p(x; t), \mu(x; t), q) = \ip\bigl[U_{\xi_*}'(y-q)m(t, y+\xi t+a)+M'_{\xi_*}(y-q)u(t, y+\xi t+a)\bigr]\,dy
\]
is of class $C^1$. Hence $\tilde q$, uniquely determined in $\mathcal V$ by the implicit function theorem, inherits the same regularity.

\end{proof}

\subsection{Analysis in the moving frame with speed \texorpdfstring{$\xi_*$}{}}

We put now ourselves in a moving frame with speed $\xi_*$, which is the speed selected by Lemma \ref{lem:xizeta}. Let $\tilde q$ be as in the previous lemma, and $q$ such that
\[
\xi_* t+q(t)=\xi t+\tilde q(t)
\]
and
\begin{align}
\theta(t,x)&=u(t,x+\xi_* t+a+q(t))-U_{\xi_*}(x), \notag \\ \label{eq:thetaeta}
 \eta(t,x)&=m(t,x+\xi_* t+a+q(t))-M_{\xi_*}(x). 
\end{align}
Note that we have
\begin{align*}
P_* & = \ip m(t,x+\xi_* t+a +q(t))u_x(t,x+\xi_* t+a+q(t))\,dx \\
& = \ip\Big(\eta(t,x) + M_{\xi_*}(x)\Big) \Big(\theta_x(t,x) + U'_{\xi_*}(x)\Big)\, dx \\
& = P_{\xi_*} + \ip M_{\xi_*}(x)\theta_x(t,x) + \eta(t,x)U'_{\xi_*}(x) +  \eta(t,x)\theta_x(t,x) \, dx ,
\end{align*}
hence
\begin{equation}\label{eq:momentum-residual}
 \ip\bigl(U_{\xi_*}'\eta+M_{\xi_*}\theta_x+\eta\theta_x\bigr)\,dx=0
\end{equation}
in view of \eqref{eq:matched-speed}. Moreover, by \eqref{eq:crucial-identity},
\begin{equation}\label{eq:quadratic-pairings}
 \ip U_{\xi_*}'\eta\,dx=\ip M_{\xi_*}\theta_x\,dx
                    =-\frac12\ip\eta\theta_x\,dx .
\end{equation}

Let us now consider the equations solved by the pertubation $(\theta, \eta)$ defined in \eqref{eq:thetaeta} :
\begin{equation}\label{eq:perturbation}
\begin{cases}
 -\theta_t-\nu\theta_{xx}+(U_{\xi_*}'+\xi_*)\theta_x
       +\frac12\theta_x^2+q'(U_{\xi_*}'+\theta_x)=K[\eta]+\bar\lambda_{\xi_*},\\
 \eta_t-\nu\eta_{xx}-\partial_x((U_{\xi_*}'+\xi_*)\eta)
       -\partial_x((M_{\xi_*}+\eta)\theta_x)
       -q' \partial_x(M_{\xi_*}+\eta)=0.
\end{cases}
\end{equation}

We now control the derivative of the correction $q$.

\begin{proposition}\label{prop:modulation} There exists $r_0, \xi_0 > 0$ such that, if $r \le r_0$ and $|\xi| \le \xi_0$, then there exists  $C > 0$ such that the function $q$ satisfies
\[
|q'(t)| \le C (\|\eta\|_{L^2} + \|\theta_x\|_{L^2})\qquad \forall t \in [0,T],
\]
where $C$ depends only on $\xi_{\max}$.
\end{proposition}

\begin{proof} Let $r \le r_0$ and $|\xi| \le \xi_0$ as in the previous lemma. We will possibly shrink $r_0$ below. Consider the quantity
\begin{multline*}
-\ip U_{\xi_*}''(M_{\xi_*}+\eta)\,dx
             +\ip M_{\xi_*}'(U_{\xi_*}'+\theta_x)\,dx = 
 2 \int M_{\xi_*}'U_{\xi_*}' -\ip U_{\xi_*}'' \eta\,dx
             +\ip M_{\xi_*}' \theta_x\,dx \le \\
             -\frac{2\nu \underline R^2}{1-\underline R^2} -\ip U_{\xi_*}'' \eta\,dx
             +\ip M_{\xi_*}' \theta_x\,dx
\end{multline*}
where the last inequality is obtained exactly as in \eqref{eq:1000}. By Lemma \ref{lem:phase}, $\|\eta\|_{L^2} + \|\theta_x\|_{L^2} \le Cr$, hence
\begin{equation*}
-\ip U_{\xi_*}'' \eta\,dx
             +\ip M_{\xi_*}' \theta_x\,dx \le C_0 \sqrt {2\pi}C r
\end{equation*}
and therefore
\begin{equation}\label{eq:1001}
-\ip U_{\xi_*}''(M_{\xi_*}+\eta)\,dx
             +\ip M_{\xi_*}'(U_{\xi_*}'+\theta_x)\,dx \le -\frac{\nu \underline R^2}{1-\underline R^2}.
\end{equation}
by taking $r$ small enough. We now use this inequality to estimate $q'$ : differentiating in time \eqref{eq:quadratic-pairings} yields
\[
 0=\ip(U_{\xi_*}'\eta_t+M_{\xi_*}'\theta_t)\,dx.
\]
The equations for $U_{\xi_*}, M_{\xi_*}$ give
\[
 \nu U_{\xi_*}'''-(U_{\xi_*}'+\xi_*)U_{\xi_*}''
       =-K[M_{\xi_*}'],\qquad
 \nu M_{\xi_*}''+(U_{\xi_*}'+\xi_*)M_{\xi_*}'
       =-M_{\xi_*}U_{\xi_*}''.
\]
since $(K[M_{\xi_*}])_x=K[M_{\xi_*}']$. Using that $\ip\eta K[M_{\xi_*}']\,dx=\ip M_{\xi_*}'K[\eta]\,dx$
multiplying the equation for $\eta$ by $U_{\xi_*}'$ and integrating yields
\begin{align*}
 \ip U_{\xi_*}'\eta_t\,dx =-\ip M_{\xi_*}'K[\eta]\,dx-\ip (M_{\xi_*} + \eta)U_{\xi_*}''\theta_x\,dx-q' \ip U_{\xi_*}''(M_{\xi_*}+\eta)\,dx.
\end{align*}
Similarly, testing the equation for $\theta$ against $M_{\xi_*}'$ gives
\begin{align*}
 \ip M_{\xi_*}'\theta_t\,dx
 &=-\ip M_{\xi_*}U_{\xi_*}''\theta_x\,dx-\ip M_{\xi_*}'K[\eta]\,dx\\
 &\quad+\frac12\ip M_{\xi_*}'\theta_x^2\,dx
       +q' \ip M_{\xi_*}'(U_{\xi_*}'+\theta_x)\,dx,
\end{align*}
since $\ip M_{\xi_*}'\,dx=0$.
Summing gives
\begin{align*}
 0= N(t) + q'(t)
 \left[-\ip U_{\xi_*}''(M_{\xi_*}+\eta)\,dx
             +\ip M_{\xi_*}'(U_{\xi_*}'+\theta_x)\,dx\right].
\end{align*}
where
\[
N(t)=-2\ip M_{\xi_*}U_{\xi_*}''\theta_x\,dx                -2\ip M_{\xi_*}'K[\eta]\,dx 
                -\ip U_{\xi_*}''\eta\theta_x\,dx
            +\frac12\ip M_{\xi_*}'\theta_x^2\,dx
\]
By estimating $N(t)$ with the $L^2$ norms of $\theta_x$ and $\eta$, and \eqref{eq:1001}, we conclude.
\end{proof}

Let $h=\eta/M_{\xi_*}$, and define the time-dependent functions
\[
\mathcal L(t) = \ip\theta\eta\,dx + \nu \|h\|_{\xi_*}^2,\qquad
\mathcal Q(t)=\Q_{\xi_*}\bigl(\theta(t),h(t)\bigr),
\]
where $\Q_{\xi_*}$ is defined in \eqref{eq:Qdef}. Our purpose is to show that $\mathcal L$ plays the role of a Lyapunov function, that is, it is strictly decreasing in time. 

Notice that $h(t)\in H^1_{\xi_*}$, since $\ip\eta\,dx=0$. On the other hand, $\Q_{\xi_*}$ and \eqref{eq:phase-constraint} 
depend on $\theta$ only through $\theta_x$, and $\mathcal L$ does not change if a constant (depending on $t$) is added to 
$\theta$, again because $\ip\eta\,dx=0$; thus we can assume that $\ip M_{\xi_*}\theta\,dx=0$, that is $\theta(t)\in H^1_{\xi_*}$, whenever needed. Moreover, \eqref{eq:quadratic-pairings} gives
\[
\ip M_{\xi_*}\bigl(U_{\xi_*}'h-\theta_x\bigr)\,dx=0,
\]
that is, \eqref{eq:phase-constraint} holds for every $t$. In 
particular, 

Theorem~\ref{thm:travel-gap}(c) yields $\mathcal Q(t)\ge c_0\,\widetilde\Psi(t)$, where the following notation is used
\begin{align*}
\Psi(t) &= \|h(t)\|^2_{\xi_*} + \|\theta_x(t)\|^2_{\xi_*}, \\
\widetilde \Psi(t) &= \|h_x(t)\|^2_{\xi_*} + \|\theta_x(t)\|^2_{\xi_*}.
\end{align*}
Before we proceed, we recall the following Poincar\'e inequality.

\begin{lemma}\label{lem:weighted} 
For $\ip M_{\xi_*} f\,dx=0$,
\begin{equation}\label{eq:poincare}
 \ip M_{\xi_*} f^2\,dx\le C_0^2\ip M_{\xi_*}|f_x|^2\,dx.
\end{equation}
Consequently, 
\begin{equation}\label{eq:comparisons}
 \Psi(t) \le C\widetilde \Psi(t),\qquad |\mathcal L(t)|\le C\Psi(t),
\end{equation}
on $[0,T]$, for some $C$ depending on $C_0$ and $\nu$.
\end{lemma}

\begin{proof}
The proof of the weighted Poincar\'e inequality is standard, we give the details for the reader's convenience. Let $\bar f=(2\pi)^{-1}\ip f\,dx$. Since $\ip M_{\xi_*} f\,dx=0$, we have $\ip M_{\xi_*}(f-\bar f)^2\,dx=\ip M_{\xi_*} f^2\,dx+\bar f^2$.
Using standard Poincar\'e inequality on the torus gives
\[
 \ip M_{\xi_*} f^2\,dx\le\ip M_{\xi_*}(f-\bar f)^2\,dx
                    \le C_0^2\ip M_{\xi_*}|f_x|^2\,dx.
\]
Applying this to $h=\eta / M_{\xi_*}$ proves the first inequality. The second one follows by Young's inequality :
\[
 \left|\ip M_{\xi_*}\theta h\,dx\right|
 \le \frac{\|h\|_{\xi_*}^2}{2} + \frac{\|\theta\|_{\xi_*}^2}{2}
 \le \frac{C_0^2}2\Psi(t).
\]
\end{proof}

We also show that $\Psi$ controls other quantities, this will be used below.

\begin{lemma}\label{lem:energy-estimates0}
There exists $r_0$ such that, if $r \le r_0$ then
\begin{equation}
 |q'|\le C\sqrt\Psi, \qquad \frac d{dt}\|h\|_{\xi_*}^2\le C\Psi, \qquad -\frac d{dt}\|\theta_x\|_{L^2}^2\le C\Psi.
 \label{eq:fixed-directional}
\end{equation}
for some $C > 0$ independent of $T$.
\end{lemma}

\begin{proof}
Regarding the inequality for $|q'|$, it is a consequence of Proposition \ref{prop:modulation}. The one on $\frac d{dt}\|h\|_{\xi_*}^2$ comes from \eqref{eq:100}, which gives
\begin{equation*}
\frac{d}{dt} \|h\|_{\xi_*}^2 + 2\nu\|h_x\|_{\xi_*}^2=  - \ip 2M_{\xi_*} \theta_x h_x - 2M'_{\xi_*} h q' + 2M_{\xi_*} h h_x \theta_x- M'_{\xi_*} h^2 q' \, dx.
\end{equation*}
The most delicate term is the third one, which can be estimated by $C\|h\|_\infty \|h_x\|_{\xi_*} \| \theta_x\|_{\xi_*}$; $\|h\|_\infty$ is controlled by $C \|h_x\|_{\xi_*}$ by Morrey inequality, so the whole term can be reabsorbed by the left-hand side provided that $\|\theta_x\|_{\xi_*}$ (and hence $r$) is small enough.

The last inequality on $-\frac d{dt}\|\theta_x\|_{L^2}^2$ follows by testing the equation for $\theta$ by $-\theta_{xx}$ :
\[
 -\frac12\frac d{dt}\|\theta_x\|_{L^2}^2+\nu\|\theta_{xx}\|_{L^2}^2
 =-\frac12\ip U_{\xi_*}''\theta_x^2\,dx
   -q'\ip U_{\xi_*}''\theta_x\,dx
   +\ip(K[\eta])_x\theta_x\,dx .
\]
Its right-hand side is bounded by $C\Psi$.
\end{proof}

\medskip
We now evaluate $\mathcal L'$.

\begin{proposition}\label{prop:energy}
Assume that $r \le r_0$ and $|\xi| \le \xi_0$, where $r_0, \xi_0 > 0$ are as in Lemma \ref{lem:phase}. Then the function $\mathcal L$ satisfies
\begin{multline}\label{eq:energy}
\mathcal L' (t) + \mathcal Q(t) = -2\nu \ip \eta\theta_x h_x\, dx -\frac12\ip\eta\theta_x^2\,dx  \\
+ q'\left(2J_{\xi_*}\ip h\,dx - 2 \nu \ip \eta h_x \, dx+\ip\eta\theta_x\,dx \right).
\end{multline}
\end{proposition}

\begin{proof}
First, $\frac{d}{dt} \ip \eta\theta\, dx=\ip(\theta_t\eta+\theta\eta_t)\,dx$, hence using the equations for $\eta,\theta$ we get
\begin{align*}
\frac{d}{dt} \ip \eta\theta\, dx&=-\|\theta_x\|_{\xi_*}^2-\frac12\ip\eta\theta_x^2\,dx-\ip K[\eta]\eta\,dx\\
 &\quad+q'\left[\ip(U_{\xi_*}'+\theta_x)\eta\,dx
                      -\ip(M_{\xi_*}+\eta)\theta_x\,dx\right].
\end{align*}
The last bracket vanishes in view of \eqref{eq:quadratic-pairings}.

Second, we evaluate $\frac{d}{dt} \|h\|_{\xi_*}^2= 2\ip \eta_t h\, dx$  (recall that $\eta=M_{\xi_*}h$). Since the equation for $M_{\xi_*}$ \eqref{eq:flux} gives 
$\nu\eta_x+(U_{\xi_*}'+\xi_*)\eta=\nu M_{\xi_*}h_x+J_{\xi_*}h$, we have
\[
 \eta_t=\partial_x\bigl(\nu M_{\xi_*}h_x+J_{\xi_*}h+
                     (M_{\xi_*}+\eta)(\theta_x+q')\bigr).
\]
Testing by $h$ and integrating by parts gives
\begin{equation}\label{eq:100}
\frac{d}{dt} \|h\|_{\xi_*}^2= -2\nu\|h_x\|_{\xi_*}^2 - 2\ip (M_{\xi_*}+\eta)(\theta_x+q') h_x\, dx,
\end{equation}
hence
\begin{multline*}
\mathcal L' = -\|\theta_x\|_{\xi_*}^2-\frac12\ip\eta\theta_x^2\,dx-\ip K[\eta]\eta\,dx -2\nu^2 \|h_x\|_{\xi_*}^2 \\ - 2\nu \ip (M_{\xi_*}+\eta)(\theta_x+q') h_x\, dx  \\
= -\frac12\ip\eta\theta_x^2\,dx -\|\theta_x\|_{\xi_*}^2-2\nu \ip M_{\xi_*}\theta_x h_x + \eta\theta_x h_x\, dx- \\ 
\ip K[\eta]\eta\,dx -2\nu^2 \|h_x\|_{\xi_*}^2 - 2\nu q' \ip (M_{\xi_*}+\eta) h_x\, dx 
\\ = -\|\theta_x + \nu h_x\|_{\xi_*}^2 - \nu^2 \|h_x\|_{\xi_*}^2 - \ip K[\eta]\eta\,dx\\ -2\nu \ip \eta\theta_x h_x\, dx  - 2\nu q' \ip (M_{\xi_*}+\eta) h_x\, dx  -\frac12\ip\eta\theta_x^2\,dx,
\end{multline*}
and the first three terms constitute $-\mathcal Q$, since $\ip K[\eta]\eta\,dx=-\kappa(C[\eta]^2+S[\eta]^2)$. 

Finally, since $q'$ is of linear order in terms of $\eta, \theta$ (hence are linear with respect to $r$ by Proposition \ref{prop:modulation}), all the remaining terms are cubic, with the exception of $q' \ip M_{\xi_*}h_x \, dx$. The last step shows that in fact this is proportional to $J_{\xi_*}$ (and hence to $\xi_*$). Recall that \eqref{eq:quadratic-pairings} gives $\ip U_{\xi_*}'\eta\,dx=\ip M_{\xi_*}\theta_x\,dx=-\frac12\ip\eta\theta_x\,dx$.
Now, multiply the equation for $M_{\xi_*}$ in \eqref{eq:flux} by $h$ and integrate to obtain (recall that $\int \eta = 0$)
\[
 \nu\ip M_{\xi_*}'h\,dx
 =J_{\xi_*}\ip h\,dx-\ip U_{\xi_*}'\eta\,dx-\xi_*\ip\eta\,dx
 =J_{\xi_*}\ip h\,dx+\frac12\ip\eta\theta_x\,dx.
\]
Plugging the identity in the previous inequality concludes the proof.
\end{proof}

\begin{lemma}\label{lem:energy-estimates}
There exists $r_0, \xi_0 > 0$ such that, if $r \le r_0$ and $|\xi| \le \xi_0$, then
\begin{align}
\mathcal L'  + \frac{c_0}2 \widetilde \Psi&\le0. \label{eq:dissipation}
\end{align}
\end{lemma}

\begin{proof} Assume that $r \le r_0$ and $|\xi| \le \xi_0$, where $r_0, \xi_0 > 0$ are as in Lemma \ref{lem:phase}. $r_0$ and $\xi_0$ might be taken smaller in what follows.

We first show that the right-hand side of \eqref{eq:energy} can be controlled by $\Psi, \widetilde \Psi$. Indeed, all the integrals can be controlled by the $L^2$ norms of $\theta_x, h, \eta$, that are bounded by $C \sqrt{\Psi}$, while the $L^2$ norm of $h_x$ is bounded by $C \sqrt{\widetilde \Psi}$. To control the first terms, we also need to bound $\eta = h M_{\xi_*}$ in $L^\infty$, which is possible by Morrey's inequality (note that $h(t, \cdot)$ vanishes somewhere on $\T$ for every $t$).
Therefore, for some $C_1 > 0$ we have
\[
\mathcal L' (t) + \mathcal Q(t) \le C_1\Big( \widetilde \Psi(t) \Psi(t)^{1/2} + \widetilde \Psi(t)^{1/2} \Psi(t)+ |q'|(|J_{\xi_*}|  \Psi^{1/2}(t) + \widetilde \Psi^{1/2}(t) \Psi^{1/2}(t) +\Psi(t)) \Big).
\]
The estimate on $q'$ from Proposition \ref{prop:modulation} and Theorem \ref{thm:travel-gap} now give
\[
\mathcal L' (t) + c_0 \widetilde \Psi(t) \le C_2 \widetilde \Psi(t)\Big(\Psi^{1/2}(t) + |\xi_*|\Big),
\]
where $|J_{\xi_*}|\le|\xi_*|/(2\pi)$ is used in view of Lemma \ref{lem:flux}. On the one hand,
\[
\Psi(t) \le C_0(\|\eta\|_{L^2}^2 + \|\theta_x\|_{L^2}^2) \le C^2_3 r^2
\]
by \eqref{eq:boh}. On the other hand,
\[
|\xi_*| \le |\xi_*-\xi| + |\xi| \le C_4 r + \xi_0
\]
by Lemma \ref{lem:xizeta}.
Choosing now $\xi_0 \le c_0/(4C_2)$ and $r \le c_0/(4C_2(C_3+C_4))$ we obtain the first desired inequality.
\end{proof}

\begin{lemma}\label{lem:fixed-weighted}
Assume that $T \ge 1$, $r \le r_0$ and $|\xi| \le \xi_0$, where $r_0, \xi_0 > 0$ are as in Lemma \ref{lem:energy-estimates}. Then,
there are $C,\omega>0$ independent of $T$, such that
\begin{align}
 \Psi(t)&\le C\left(\|h(0)\|_{\xi_*}^2e^{-\omega t}
                 +\|\theta_x(T)\|_{\xi_*}^2e^{-\omega(T-t)}\right),
                 \label{eq:fixed-pointwise}\\
\int_0^T|q'(t)|\,dt
 &\le C\left(\|h(0)\|_{\xi_*}+\|\theta_x(T)\|_{\xi_*}\right).
                 \label{eq:fixed-integrated-phase}
\end{align}
\end{lemma}
\begin{proof}

Let us first notice that \eqref{eq:comparisons} gives the existence of $C_1$ such that $|\mathcal L(t)| \le C_1 \widetilde \Psi(t)$. Choose $\omega \le c_0/(4C_1)$. For any positive Lipschitz $w$ with $|w'|\le\omega w$,
integrating \eqref{eq:dissipation} times $w$ gives
\begin{equation}\label{eq:200}
\frac{c_0}2 \int_0^T w \widetilde \Psi \, dt \le w(0)\mathcal L(0)-w(T)\mathcal L(T) + \int_0^T w' \mathcal L \, dt \le w(0)\mathcal L(0)-w(T)\mathcal L(T) + \frac{c_0}4 \int_0^T w \widetilde \Psi \, dt.
\end{equation}
For any $\epsilon>0$, Poincar\'e and Young inequalities give
\begin{align*}
 \mathcal L(0)\le & \epsilon \|\theta_x(0)\|^2_{L^2} + C_\epsilon\|h(0)\|^2_{L^2}\\
 -\mathcal L(T)\le & -\frac{\nu} {2C_0}\|h(T)\|^2_{L^2}+C_2\|\theta_x(T)\|^2_{L^2}.
\end{align*}
We now integrate the last inequality in \eqref{eq:fixed-directional}, first on $[0,t]$ then on $[0,1]$ to get
\[
\|\theta_x(0)\|^2_{L^2} \le \int_0^1\|\theta_x(t)\|^2_{L^2} \, dt + C_3\int_0^1 \widetilde \Psi \, dt.
\]
Therefore,
\[
w(0) \|\theta_x(0)\|^2_{L^2} \le (C_3 + 1)e^\omega\int_0^1 w \widetilde \Psi \, dt.
\]
Plugging everything into \eqref{eq:200} gives
\[
\frac{c_0}4 \int_0^T w \widetilde \Psi \, dt \le \epsilon(C_3 + 1)e^\omega\int_0^1 w \widetilde \Psi \, dt + w(0) C_\epsilon \|h(0)\|^2_{L^2} + w(T) C_2\|\theta_x(T)\|^2_{L^2},
\]
that is, for $\epsilon$ small enough,
\[
\int_0^T w \widetilde \Psi \, dt \le C_3[w(0) \|h(0)\|^2_{L^2} + w(T) \|\theta_x(T)\|^2_{L^2}].
\]
For $w(t)=e^{-\omega|t-\tau|}$ this controls the integral of $\widetilde\Psi$ near $\tau$ (say, on $[(\tau-1/2) \vee 0, (\tau+1/2) \wedge T]$). Integrating \eqref{eq:fixed-directional} give the pointwise estimate on $\|h(\tau)\|_{\xi_*}^2$ and $\|\theta_x(\tau)\|_{\xi_*}^2$, and therefore on $\Psi(\tau)$, that is \eqref{eq:fixed-pointwise}. The following inequality \eqref{eq:fixed-integrated-phase} follows by further integration on $[0,T]$ of $\sqrt \Psi$ and the first inequality in \eqref{eq:fixed-directional}.
\end{proof}

\begin{theorem}\label{thm:half} Let $r_0, \xi_0 > 0$ be as in Lemma \ref{lem:energy-estimates}. Then, there exists $e_0, T_0 > 0$ such that, if $T \ge T_0$,
\begin{align*}
\sup_{t \in [0,T]} \left(\|m(t, \cdot)-M_{\xi}(\cdot-\xi t - a)\|_{L^2}+\|u_x(t, \cdot)-U'_{\xi}(\cdot-\xi t - a)\|_{L^2}\right) & \le r_0,\\
\|m_0-M_{\xi}(\cdot-a)\|_{L^2} +\|\partial_x u_T-U'_{\xi}(\cdot-\xi T-a)\|_{L^2} & \le e_0,
\end{align*}
for some $a$ and $|\xi| \le \xi_0$, then 
\begin{equation}\label{eq:improvement}
\sup_{t \in [0,T]} \left(\|m(t, \cdot)-M_{\xi}(\cdot-\xi t - a)\|_{L^2}+\|u_x(t, \cdot)-U'_{\xi}(\cdot-\xi t - a)\|_{L^2}\right) \le \frac {r_0} 2.
\end{equation}
Moreover, there exists $\xi_*$ and $q \in C^1([0,T])$ with the property that
\[
|\xi_*-\xi| \le \frac{C}{T}e_0.
\]
and
\begin{multline}\label{eq:turnpikee}
|q'(t)|+\|m(t, \cdot)-M_{\xi_*}(\cdot-\xi_* t - a-q(t))\|_{L^2}+\|u_x(t, \cdot)-U'_{\xi_*}(\cdot-\xi_* t - a-q(t))\|_{L^2} \le \\ Ce_0\left(e^{-\omega t}+e^{-\omega(T-t)}\right) \qquad \forall t \in [0,T]
\end{multline}
for some $C, \omega > 0$.

\end{theorem}

\begin{proof} Recall that
\begin{align*}
\theta_x(t,x)&=u_x(t,x+\xi_* t+a+q(t))-U'_{\xi_*}(x) = u_x(t,x+\xi t+a+\tilde q(t))-U'_{\xi_*}(x)\\
 \eta(t,x)&=m(t,x+\xi_* t+a+q(t))-M_{\xi_*}(x) = m(t,x+\xi t+a+\tilde q(t))-M_{\xi_*}(x)
\end{align*}
where $\xi_*, q, \tilde q$ are given by the above Lemmata. Recall that $|\xi_*- \xi| \le C r_0$ by Lemma \ref{lem:xizeta}. A crucial step will be to show that in fact $\xi_*- \xi$ is of order $O(1/T)$. Let us first observe that by \eqref{eq:phase-residual}
\[
|\tilde q(0)| \le C_1\left|\ip\bigl[U_{\xi_*}'(x)(m(0, x+a)-M_{\xi_*}(x))-M_{\xi_*}(x)(u_x(0, x+a)-U_{\xi_*}'(x))\bigr]\,dx \right|.
\]
To estimate the right-hand side with $e_0$, argue as in \eqref{eq:123} but starting from $0=P_*-P_{\xi_*}=...$ to get
\begin{multline*}
-\ip (u_x(0, x+a) - U'_{\xi_*}) M_{\xi_*}\, dx= \\ \ip (m(0, x+a)-M_{\xi_*})(u_x(0, x+a) - U'_{\xi_*}) + (m(0, x+a)-M_{\xi_*})U'_{\xi_*} \, dx,
\end{multline*}
so
\[
|\tilde q(0)| \le C_1\|m(0, \cdot+a)-M_{\xi_*}\|_{L^2}\|u_x(0, \cdot+a) + U'_{\xi_*}\|_{L^2}\le C_2(e_0 + |\xi_*-\xi|)
\]
using \eqref{eq:noescape-translate-profiles}. An analogous estimate holds for $|\tilde q(T)|$.

Now we want to control $q$ using \eqref{eq:fixed-integrated-phase}. Again using \eqref{eq:noescape-translate-profiles},
\[
\|h(0)\|_{\xi_*} \le C_0^{1/2}\|m(0,\cdot+a+\tilde q(0))-M_{\xi_*}(x)\|_{L^2} \le C_3(e_0 + |\xi_*-\xi|),
\]
and an analogous estimate holds for $\|\theta_x(T)\|_{\xi_*}$. By Lemma \ref{lem:fixed-weighted} we get
\begin{align}
 \Psi(t)&\le C_4(e_0 + |\xi_*-\xi|)^2\left(e^{-\omega t}
                 +e^{-\omega(T-t)}\right),
                 \label{eq:fixed-pointwiseb}\\
\int_0^T|q'(t)|\,dt
 &\le C_4(e_0 + |\xi_*-\xi|).
                 \label{eq:fixed-integrated-phaseb}
\end{align}
Since $t(\xi_*-\xi) = \tilde q(t)-q(t)$ for all $t$,
\[
T|\xi_*-\xi| \le |\tilde q(T) - \tilde q(0)| + |q(T) - q(0)| \le (2C_2+C_4)(e_0 + |\xi_*-\xi|).
\]
Whenever $T \ge 2(2C_2+C_4)$, we have
\[
|\xi_*-\xi| \le \frac{2(2C_2+C_4)}{T}e_0 \quad ( \le e_0).
\]
Plugging this inequality into the previous ones yields
\[
|\tilde q(0)| \le 2C_2 e_0, \quad \sup_t \sqrt \Psi(t) \le 2\sqrt {2 C_4} e_0 \quad \int_0^T|q'(t)|\,dt \le 2C_4 e_0
\]
so, by the fundamental theorem of calculus,
\[
|\tilde q(t)| \le |\tilde q(0)| + t|\xi_*-\xi| + |q(t) - q(0)| \le C_5 e_0.
\]
Finally, by \eqref{eq:noescape-translate-profiles}
\begin{multline*}
\|m(t, \cdot)-M_{\xi}(\cdot-\xi t - a)\|_{L^2} = \|m(t, \cdot+\xi t + a + \tilde q)-M_{\xi}(\cdot + \tilde q)\|_{L^2} \\ \le
\|\eta\|_{L^2} + \|M_{\xi_*}-M_{\xi}(\cdot + \tilde q)\|_{L^2} \le C_6(\sqrt \Psi + |\tilde q(t)| + |\xi_*-\xi|) \le C_7 e_0
\end{multline*}
for all $t \in [0,T]$. An analogous bound can be obtained for $\|u_x(t, \cdot)-U'_{\xi}(\cdot-\xi t - a)\|_{L^2}$. Therefore, choosing $e_0 \le r_0/(4C_7)$ gives \eqref{eq:improvement}. The inequality \eqref{eq:turnpikee} is a straightforward consequence of \eqref{eq:fixed-pointwiseb} and \eqref{eq:fixed-directional}.
\end{proof}

\begin{remark}\label{rem:unique} In the conclusions of the previous theorem one can in fact replace $e_0$ by $\|m_0-M_{\xi}(\cdot-a)\|_{L^2} +\|\partial_x u_T-U'_{\xi}(\cdot-\xi T-a)\|_{L^2}$. As a straightforward consequence, if $m_0=M_{\xi}(\cdot-a)$ and $\partial_x u_T=U'_{\xi}(\cdot-\xi T-a)$, then $\xi_* = \xi$ and $m(t, \cdot)=M_{\xi}(\cdot-\xi t - a)$ and $u_x(t, \cdot)=U'_{\xi}(\cdot-\xi t - a)$ for all $t$. In other words, in a neighborhood of the traveling wave, the traveling wave itself is the unique solution with initial/final conditions $M_{\xi}(\cdot-a), U'_{\xi}(\cdot-\xi T-a)$.
\end{remark}

\begin{remark}\label{rem:momentum} It is clear that the a priori the selected speed $\xi_*$ might be different from the reference one $\xi$. In fact, they coincide provided that
\[
 P_*=\ip m(t,x)u_x(t,x)\,dx = P_{\xi},
\]
in which case there is no need to use Lemma \ref{lem:xizeta}. This equality may come from further assumptions. For instance, if $u_T = 0$, then $P_* = 0$, which naturally forces comparison with the stationary profiles $(U_0, M_0)$. Similarly, assume that we are in the even setting, that is we work with solutions $(u,m)$ that are symmetric around $x = 0$ (or around any other $x = a$) Since $P_* = 0$, we get the standard turnpike property
\[
\|m(t, \cdot)-M_0\|_{L^2}+\|u_x(t, \cdot)-U'_0\|_{L^2} \le 
 Ce_0\left(e^{-\omega t}+e^{-\omega(T-t)}\right)
\]
since there is no need to implement the additional modulation $q$.
\end{remark}

\begin{remark}\label{rem:uniquexi} The selected speed $\xi_*$ of Theorem \ref{thm:half} is uniquely determined by the momentum argument of Lemma \ref{lem:xizeta}. However, it is not the unique one for which the exponential estimate \eqref{eq:turnpikee} holds. Indeed, consider $\bar \xi$ close to $\xi_*$. We have, setting $\bar q = q + (\xi_*- \bar \xi)t$ (and dropping the comparison between $u_x$ and $U'_{\bar \xi}$ for simplicity)
\begin{multline*}
|\bar q'(t)|+\|m(t, \cdot)-M_{\bar \xi}(\cdot-\bar \xi t - a-\bar q(t))\|_{L^2} \le \\
|q'(t)|+\|m(t, \cdot)-M_{\bar \xi}(\cdot-\xi_* t - a-q(t))\|_{L^2} + |\xi_* - \bar \xi | \le \\
|q'(t)|+\|m(t, \cdot)-M_{\xi_*}(\cdot-\xi_* t - a-q(t))\|_{L^2} + (C_0+1) |\xi_* - \bar \xi | \le
\\ (C_0+1) |\xi_* - \bar \xi |+ Ce_0\left(e^{-\omega t}+e^{-\omega(T-t)}\right),
\end{multline*}
thus
\[
|\bar q'(t)|+\|m(t, \cdot)-M_{\bar \xi}(\cdot-\bar \xi t - a-\bar q(t))\|_{L^2} \le C_1e_0\left(e^{-\omega t}+e^{-\omega(T-t)}\right)
\]
holds provided that 
\[
|\xi_* - \bar \xi | \lesssim e_0\,e^{-\omega T /2}.
\]

\end{remark}

\subsection{Invertibility of the linearized operator}

To conclude this section, we show that the linearization of \eqref{eq:mfg} around the wave 
$U_{\xi}(x - \xi t), M_{\xi}(x - \xi t)$ 
\begin{equation}
 \begin{cases}
 -p_t-\nu p_{xx}+(U_{\xi}'+\xi)p_x=K[n],&p(T)=0,\\
 n_t-\nu n_{xx}-\partial_x((U_{\xi}'+\xi)n)
                 =\partial_x(M_{\xi} p_x),&n(0)=0 .
 \end{cases}
 \label{eq:continuation-linearization}
\end{equation}
has only the trivial solution. This will be used in the existence theorem of the next section. The 
arguments follow the lines of the ones used to get to Theorem \ref{thm:half}, tailored to the linearized 
framework. In particular, the linearized momentum
\[
\ip(M_{\xi} p_x+U_{\xi}'n)\,dx
\]
will be used. The main difference with respect to the nonlinear framework will be that, instead of using 
the change of variables $x \mapsto x + \tilde q(t)$ to align the momentum of the dynamic solution with 
the one of a traveling wave with (slightly) different speed $\xi_*$, an additive correction will be 
used. This is actually coherent with the first-order expansion
\[
(M_{\xi}+\epsilon n)(t,\cdot +\epsilon \beta(t)) - M_{\xi + \epsilon \gamma} = \epsilon[n(t,\cdot +\epsilon \beta(t)) + \beta M'_{\xi} - \gamma \partial_{\xi} M_{\xi}] + o(\epsilon).
\]

\begin{proposition}\label{prop:linvanish} Let $\xi_0$ be as in Theorem \ref{thm:half}, and $|\xi|\le\xi_0$. Then, the only classical solution to \eqref{eq:continuation-linearization} is $(n, p) = (0,0)$, provided that $T \ge T_0$ (which can be larger than the one appearing in Theorem \ref{thm:half}).
\end{proposition}

\begin{proof}

First, differentiating \eqref{eq:continuation-linearization} as in Proposition~\ref{prop:momentum} gives a constant (linearized) momentum
\[
 \ell:=\ip(M_{\xi} p_x+U_{\xi} 'n)\,dx .
\]
For $\beta=\beta(t)$ and $\gamma \in \R$, let
\[
 \begin{aligned}
 \eta&:=n-\gamma\partial_{\xi} M_{\xi}+\beta M_{\xi}',\\
 \theta&:=p-\gamma \partial_{\xi} U_{\xi}+\beta U'_{\xi}
 \end{aligned}
\]
Choosing $\gamma:=\ell/P_{\xi}'$ (which is well defined in view of \eqref{eq:working})  we have
\[
 \ip(U_{\xi}'\eta+M_{\xi}\theta_x)\,dx=\ell-\gamma P_{\xi}'=0.
\]
On the other hand,
\[
\ip(U_{\xi}'\eta -M_{\xi}\theta_x)\,dx=2\ip U'_{\xi}\eta\,dx=2\ip U'_{\xi}(n-\gamma \partial_{\xi} M_{\xi})\,dx + 2\beta \ip M'_{\xi}  U'_{\xi}\,dx = 0
\]
whenever
\[
\beta(t) = -\frac{\ip U_{\xi}'(n-\gamma \partial_{\xi} M_{\xi})\,dx}{\ip M_{\xi}'  U_{\xi}'\,dx} \quad \left( \, = -\frac{\ip M_{\xi}'(p-\gamma \partial_{\xi} U_{\xi})\,dx}{\ip M_{\xi}'  U_{\xi}'\,dx}\, \right),
\]
which is well defined by \eqref{eq:1000}. By this choice of $\gamma, \beta$, both $\int U_{\xi}'\eta$ and $\int M_{\xi}\theta_x$ vanish.

Differentiating the system \eqref{eq:profiles} for $U_{\xi}, M_{\xi}$ with respect to $\xi$ and $x$ shows that
\[
 \begin{cases}
 -\theta_t-\nu\theta_{xx}+(U_{\xi}'+\xi)\theta_x+(\beta'-\gamma) U_{\xi}'
                         =K[\eta] +\gamma\,\partial_{\xi}\bar\lambda_{\xi},\\
 \eta_t-\nu\eta_{xx}-\partial_x((U_{\xi}'+\xi)\eta)
                 -\partial_x(M_{\xi}\theta_x)-(\beta'-\gamma) M_{\xi}'=0 .
 \end{cases}
\]
Differentiating in time the identity $\ip(U_{\xi}'\eta +M'_{\xi}\theta)\,dx=0$, as in
Proposition~\ref{prop:modulation}, gives 
\[
2(\beta'-\gamma)\ip M'_{\xi} U'_{\xi} \, dx
   =2\ip M_{\xi} U_{\xi}''\theta_x\,dx+2\ip M_{\xi}'K[\eta]\,dx,
\]
hence, similarly as above, for $h=\eta/M_{\xi}$ and $\Psi=\|h\|_{\xi}^2+\|\theta_x\|_{\xi}^2$ and
$\widetilde\Psi=\|h_x\|_{\xi}^2+\|\theta_x\|_{\xi}^2$ we get
\[|\beta'-\gamma|\le C\sqrt\Psi . \]

Arguing as in Proposition \ref{prop:energy}, that is, testing the equations for $\eta, \theta$  gives the identity
\[
 \frac{d}{dt}\left(\ip\theta\eta\,dx+\nu\|h\|_{\xi}^2\right)
       + \Q_{\xi}(\theta,h) =2J_{\xi} (\beta'-\gamma) \ip h\,dx .
\]
Theorem~\ref{thm:travel-gap}, Lemma~\ref{lem:weighted}, and
$|J_{\xi}|\le|\xi|/(2\pi)$ imply 
\[
 \mathcal L'+\frac {c_0} 2\widetilde\Psi\le0,\qquad
 \frac d{dt}\|h\|_{\xi}^2\le C\Psi,\qquad
 -\frac d{dt}\|\theta_x\|_{L^2}^2\le C\Psi.
\]
The last two inequalities follow by the same steps
of Lemma~\ref{lem:energy-estimates0}, with no nonlinear terms: these are quadratic estimates for a linear system, so no smallness of solutions is required.
The proof of Lemma~\ref{lem:fixed-weighted} uses only these inequalities, weighted Poincar\'e, and $|\mathcal L|\le C\Psi$.
Repeating that proof therefore gives
\begin{equation}
 \sup_t\sqrt\Psi+\int_0^T|\beta'-\gamma|\,dt
       \le C\bigl(\|h(0)\|_{\xi}+\|\theta_x(T)\|_{\xi}\bigr),
 \label{eq:continuation-linear-bound}
\end{equation}
with $C$ independent of $T$.

We now conclude as in Theorem \ref{thm:half}. Using the definition of $\beta$, and $n(0)=0$, $p(T)=0$
\[
 \beta(0)=\gamma  \frac{\ip U_{\xi}' \partial_{\xi} M_{\xi}\,dx}{\ip M_{\xi}'  U_{\xi}'\,dx},\qquad
 \beta(T)=\gamma  \frac{\ip M_{\xi}' \partial_{\xi} U_{\xi}\,dx}{\ip M_{\xi}'  U_{\xi}'\,dx},
\]
which, by \eqref{eq:noescape-translate-profiles}, give
\[
 |\beta(0)|+|\beta(T)|
       +\|h(0)\|_{\xi}+\|\theta_x(T)\|_{\xi}\le C|\gamma|.
\]
Equation \eqref{eq:continuation-linear-bound} implies
\[
 T|\gamma|
 \le|\beta(T)|+|\beta(0)|+\int_0^T|\beta'-\gamma|\,dt
 \le C|\gamma|.
\]
Taking $T_0 \ge 2C$, we obtain $\gamma=0$, and \eqref{eq:continuation-linear-bound} gives $\Psi=0$.
Then the formula for $\beta(0)$ gives $\beta(0) = 0$, which together with $|\beta'-\gamma| = 0$ implies $\beta=0$.
Thus $n=0$ and $p_x=0$. The equation for $p$ and $p(T)=0$ then give $p=0$.
\end{proof}

\subsection{Existence by continuation}
\label{sec:fixed-speed-existence}

The existence of a solution to \eqref{eq:mfg} which stays close to a traveling wave will be obtained by a 
continuation argument, based on the Leray--Schauder degree, a well-known extension of finite dimensional 
Brouwer degree. For the reader's convenience, we recall the 
facts that we need, referring to e.g. \cite[Chapter~4]{AmbroArcoBook} for further details. 
Let $X$ be a real Banach space, $\Omega\subset X$ a bounded open set, and 
$\mathcal T\colon\overline\Omega\to X$ a (nonlinear) compact map. If 
$u\neq\mathcal T(u)$ for every $u\in\partial\Omega$, the Leray--Schauder degree $\deg(I-\mathcal 
T,\Omega,0)\in\mathbb Z$ is well defined. 
\begin{theorem}[{\cite[Chapter~4]{AmbroArcoBook}}]\label{thm:LS}
Let $X$ and $\Omega$ be as above.
\begin{enumerate}
\item[(i)] Let $\mathcal T\colon[0,1]\times\overline\Omega\to X$ be a compact map, and assume that 
$u\neq\mathcal T(\tau,u)$ for every $\tau\in[0,1]$ and $u\in\partial\Omega$. Then 
$\deg(I-\mathcal T(\tau,\cdot),\Omega,0)$ does not depend on $\tau$. If, in addition, 
$\deg(I-\mathcal T(0,\cdot),\Omega,0)\neq0$, then $\mathcal T(1,\cdot)$ has a fixed point in $\Omega$.
\item[(ii)] Let $\mathcal T\colon\overline\Omega\to X$ be a compact map of class $C^1$, and let 
$u_0\in\Omega$ be its unique fixed point in $\overline\Omega$. If $I-\mathcal T'(u_0)$ is injective, then
\[
 \deg(I-\mathcal T,\Omega,0)=(-1)^\beta\neq0,
\]
where $\beta$ is the sum of the algebraic multiplicities of the eigenvalues of $\mathcal T'(u_0)$ which 
are larger than $1$.
\end{enumerate}
\end{theorem}

We will apply (i) in the theorem above to a suitable operator homotopy, for which fixed points of 
$\mathcal T(1,\cdot)$ are solutions of \eqref{eq:mfg}, which existence will follow by the nontrivial 
degree of $I-\mathcal T(0,\cdot)$. For the latter property, triggered by (ii), we stress that we will not 
need the actual value of $\beta$, but only the fact that the degree does not vanish, i.e. that the 
unique fixed point is nondegenerate. Notice that, differently from the Schauder fixed point theorem 
(which is itself a consequence of (i), see \cite[Theorem~4.2.6]{AmbroArcoBook}), no set which is 
invariant under $\mathcal T(1,\cdot)$ is needed: it is enough to know that the fixed points of 
$\mathcal T(\tau,\cdot)$ do not reach $\partial\Omega$, which is the content of an a priori estimate, 
together with the nondegeneracy of isolated fixed points of $\mathcal T(0,\cdot)$. In 
our case, the former a priori estimate is provided by Theorem~\ref{thm:half}, while the latter nondegeneracy 
comes from Proposition~\ref{prop:linvanish}.

\begin{theorem}\label{thm:continuation} Assume that $m_0,u_T\in C^{2,\alpha}(\T)$ for $\alpha \in (0,1)$, 
$m_0>0$, $\ip m_0\,dx=1$. 

Let $e_0, \xi_0, T_0$ be as in Proposition \ref{prop:linvanish}. Assume that for some $a$ and 
$|\xi| \le \xi_0$ it holds
\[
\|m_0-M_{\xi}(\cdot-a)\|_{L^2} +\|\partial_x u_T-U'_{\xi}(\cdot-\xi T-a)\|_{L^2} \le e_0.
\]
Then, for any $T \ge T_0$, \eqref{eq:mfg} has at least one classical solution satisfying 
\eqref{eq:turnpikee}.
\end{theorem}

\begin{proof}
Fix $T\ge T_0$. We will apply Theorem~\ref{thm:LS} in the space
\[
 X_T:=C([0,T];L_0^2(\T)\times L_0^2(\T)),\qquad
 \|(\rho,v)\|_{X_T}:=\sup_{t \in [0,T]}(\|\rho(t)\|_{L^2}+\|v(t)\|_{L^2}),
\]
where $L_0^2(\T)$ denotes the functions with zero integral, and with
$\Omega=B:=\{Z\in X_T:\|Z\|_{X_T}<r_0\}$. The parameter $\tau\in[0,1]$ moves the initial and final 
data from the ones of the wave $(U_{\xi}(x - \xi t-a), M_{\xi}(x - \xi t-a))$, at 
$\tau=0$, to $m_0$ and $u_T$, at $\tau=1$; it is convenient to work in the frame moving with speed 
$\xi$ and centered at $a$. 

For $\tau\in[0,1]$ and $Z=(\rho,v)\in X_T$, we first solve the backward equation and then the forward 
equation
\begin{equation}
 \begin{cases}
 -w_t-\nu w_{xx}+\xi w_x+\frac12w_x^2=K[M_{\xi}+\rho],
       &w(T)=\tau u_T(\cdot+\xi T+a) +(1-\tau) U_{\xi},\\
 \mu_t-\nu\mu_{xx}-\partial_x((\xi+w_x)\mu)=0,
       &\mu(0)=\tau m_0(\cdot+a)+(1-\tau)M_{\xi} ,
 \end{cases}
 \label{eq:continuation-response}
\end{equation}
and we define
\begin{equation}
 \mathcal H_\tau(Z):=(\mu-M_{\xi},w_x-U_{\xi}').
 \label{eq:continuation-map}
\end{equation}
Note that the density at $t=0$ is positive and it has unit mass, being a convex combination of positive 
probability densities; then $\mu(t)$ has unit mass for every $t$, and $\mathcal H_\tau$ maps $X_T$ into 
itself. Moreover, it is standard to show that the operator $\mathcal H_\tau$ is well-defined (the first 
equation in particular can be transformed into a linear one by the Hopf--Cole transformation). Since 
$K$ is a nonlocal, smoothing operator, whenever $\rho$ varies in a bounded set of $L^2$, linear 
parabolic estimates give bounded image in $C^{\alpha/2, \alpha} \times C^{\alpha/2, \alpha}$, and the 
Ascoli--Arzel\`a theorem makes $(\tau,Z)\mapsto\mathcal H_\tau(Z)$ a compact map on 
$[0,1]\times\overline B$. Its continuity is standard by considering the difference between the 
equations. Similarly, since $\rho\mapsto K[\rho]$ is linear and continuous from $L^2(\T)$ into 
$C^k(\T)$, for every $k$, and the solutions of \eqref{eq:continuation-response} depend smoothly on the 
right-hand side of the first equation (again by the Hopf--Cole transformation), the map 
$Z\mapsto\mathcal H_0(Z)$ is of class $C^1$, and its derivative is obtained by linearizing 
\eqref{eq:continuation-response}.

\smallskip
We first check that the fixed points of $\mathcal H_\tau$ do not reach the boundary of $B$. If 
$Z=\mathcal H_\tau(Z)$ and $Z\in\overline B$, then
\[
 u(t,x)=w(t,x-\xi t-a),\qquad m(t,x)=\mu(t,x-\xi t-a)
\]
solve \eqref{eq:mfg} with $m(0)=(1-\tau)M_{\xi}(\cdot-a)+\tau m_0$, 
$u(T)=(1-\tau)U_{\xi}(\cdot-a-\xi T)+\tau u_T$. Since \eqref{eq:apriori} is exactly 
$\|Z\|_{X_T} \le r_0$, and
\[
\|m(0)-M_{\xi}(\cdot-a)\|_{L^2} +\|\partial_x u(T)-U'_{\xi}(\cdot-\xi T-a)\|_{L^2} \le \tau e_0,
\]
Theorem~\ref{thm:half} applies, yielding
\begin{equation}
 Z=\mathcal H_\tau(Z),\quad \|Z\|_{X_T}\le r_0
       \quad\Longrightarrow\quad \|Z\|_{X_T}\le r_0/2 .
 \label{eq:continuation-boundary}
\end{equation}
Thus $Z\neq\mathcal H_\tau(Z)$ for every $\tau\in[0,1]$ and $Z\in\partial B$.

\smallskip
Next, we show that $\deg(I-\mathcal H_0,B,0)\ne0$. Clearly $\mathcal H_0(0)=0$, and the corresponding 
solution is $w=U_{\xi}+\bar\lambda_{\xi}(T-t)$, $\mu=M_{\xi}$, that is the wave itself. By 
\eqref{eq:continuation-boundary} and Remark~\ref{rem:unique}, $Z=0$ is the unique fixed point of 
$\mathcal H_0$ in $\overline B$. Moreover, $\mathcal H_0'(0)Z=(n,p_x)$, where $(n,p)$ solves
\[
 \begin{cases}
 -p_t-\nu p_{xx}+(U_{\xi}'+\xi)p_x=K[\rho],&p(T)=0,\\
 n_t-\nu n_{xx}-\partial_x((U_{\xi}'+\xi)n)
                 =\partial_x(M_{\xi} p_x),&n(0)=0 .
 \end{cases}
\]
Hence, if $Z=(\rho,v)$ satisfies $\mathcal H_0'(0)Z=Z$, then $\rho=n$ and $v=p_x$, and $(n,p)$ is a 
solution of the linearized system \eqref{eq:continuation-linearization}, which is classical by parabolic 
regularity. By Proposition~\ref{prop:linvanish}, $(n,p)=(0,0)$, that is $Z=0$. Therefore 
$I-\mathcal H_0'(0)$ is injective, and Theorem~\ref{thm:LS}(ii) gives $\deg(I-\mathcal H_0,B,0)\neq0$.

\smallskip
Finally, by Theorem~\ref{thm:LS}(i), $\mathcal H_1$ has a fixed point $Z\in B$. At a fixed point, 
$K[M_{\xi}+\rho]=K[\mu]$; then $(u,m)$, defined as above with $\tau=1$, solves \eqref{eq:mfg} with 
$m(0)=m_0$ and $u(T)=u_T$, and parabolic regularity shows that it is a classical solution. Since 
$\|Z\|_{X_T}<r_0$, Theorem~\ref{thm:half} applies, and gives \eqref{eq:turnpikee}.
\end{proof}

\begin{remark}[On periodic solutions]\label{rem:periodic}
Since the family of traveling waves $U_{\xi}(x-\xi t), M_{\xi}(x-\xi t)$ solving \eqref{eq:mfg} is periodic in the $x$-variable, each solution can be thought as a \emph{periodic in time} solution with period
\[
\tau_{\xi} = \frac{2\pi}{|\xi|}
\]
The previous existence result says that these periodic solutions are somewhat stable, at least when $|\xi| \le \xi_0$. Consider indeed small perturbations $f, g$, and, for every $k \in \mathbb N$ large, apply Theorem \ref{thm:continuation} to
\[
m_0 = M_{\xi} + f, \qquad u_T = U_{\xi} + g, \qquad T = T_k := k \tau_{\xi} = \frac{2\pi k}{|\xi|}.
\]
Whenever $f$, $g$ are small enough (independently on $k$), the Theorem gives the existence of 
$\xi_{*,k}$ and $q_k \in C^1([0,T])$ with the property that
\[
\frac{|\xi_{*,k}-\xi|}{|\xi|} \le \frac{C}{k}.
\]
and
\begin{multline*}
|q_k'(t)|+\|m_k(t, \cdot)-M_{\xi_{*,k}}(\cdot-\xi_{*,k} t-q_k(t))\|_{L^2}+\|\partial_x u_k(t, \cdot)-U'_{\xi_{*,k}}(\cdot-\xi_{*,k} t - q_k(t))\|_{L^2} \le \\ C\left(e^{-\omega t}+e^{-\omega(T-t)}\right).
\end{multline*}
The two inequalities say that $(u_k, m_k)$ must be exponentially close to a periodic solution with period $\tau_k$, which approaches $\tau_{\xi}$ as $k \to \infty$. It is interesting to observe what happens in the limit $k \to \infty$, assuming the validity of compactness to perform the following operations (which is expected, see for instance \cite[Theorem 4.3]{CD26}): $u_k \to \bar u$, $m_k \to \bar m$, where $(\bar u, \bar m)$ solve the MFG system on $(0,\infty) \times \mathbb T$, $\bar m(0) = M_{\xi} + f$ and $q_k \to \bar q$. One would have in particular the following long-time behavior:
\[
\|\bar m(t, \cdot)-M_{\xi}(\cdot-\xi t- \bar q)\|_{L^2}\le C e^{-\omega t}, 
\]
with $|\bar q(t) - q^*| \le C e^{-\omega t}$ for some $q^* \in \T$ (since $\bar q'(t)$ would satisfy the previous inequality). Therefore, \eqref{eq:noescape-translate-profiles} would give
\[
\|\bar m(t, \cdot)-M_{\xi}(\cdot-\xi t- q^*)\|_{L^2}\le C e^{-\omega t}.
\]

This argument suggests also the following fact on stationary solutions $(U_0,M_0)$. If we take $\bar m(0)$ as a perturbation of $M_0$, there exists a solution $(\bar u, \bar m)$ to \eqref{eq:mfg} on $(0,\infty) \times \mathbb T$ that converges exponentially fast to a \emph{translation} of $M_0$.

\medskip

There is a slightly different viewpoint on this construction; take $m_0=M_0$ itself as the initial condition and $u_T = U_0$, $T$ large enough, $N$ be any integer such that $|N|\le \xi_0T/(2 \pi)$, and $\xi_{N,T} = 2\pi N/T$. Since $\|m_0-M_\xi\|_{L^2} \le C|\xi|$, and the same holds for $u_T - U_{\xi_{N,T}}(\cdot - \xi_{N,T} T) = u_T - U_{\xi_{N,T}}$, the existence theorem applies, giving solutions exponentially close to a wave with a speed comparable to $\xi_{N,T}$. These solutions can be seen as cyclic motions, in particular they complete $N$ revolutions during the time interval $[0,T]$; $N$ plays the role of the \emph{winding number} of solutions. This shows that various possible equilibria arise as the time horizon increases.

In the $T \to \infty$ limit, these arguments also suggests the existence of a family of solutions $(m_\xi, u_\xi)$ to \eqref{eq:mfg}, with $|\xi|$ small, such that $m_\xi(0) = M_0$ for all $\xi$ and 
\[
\|m_\xi(t, \cdot)-M_{\xi}(\cdot-\xi t- q_\xi)\|_{L^2}\le C e^{-\omega t}. 
\]
That is, even if we start the evolution from the stationary state $M_0$, infinitely many equilibria traveling at small speed form.

\end{remark}

\appendix

\section{Proof of Proposition~\ref{prop:concavity}}\label{sec:concavity}

In this appendix we prove Proposition~\ref{prop:concavity}, following the strategy of
\cite[Lemmas~3.1, 3.2, 4.5 and Proposition~4.6]{munoz2026uniqueness}, with some simplifications
(there the discounted problem is considered, and our parameter $\tau$ is denoted by $\zeta$).

Throughout the appendix $\tau>0$ is fixed, and we drop the index $\tau$, writing $v$, $m$, $g$, $A$ for
$v_\tau$, $m_\tau$, $g_\tau$, $A(\tau)$. Moreover, we set
\[
 \phi:=v',\qquad \qquad Z:=\frac{g'}{\sin x}.
\]
By Lemmas~\ref{lem:scalar} and \ref{lem:corr-g}, $v$, $m$ and $g$ are even and
$2\pi$-periodic, hence they are even also with respect to $x=\pi$. As a consequence, $\phi$ and $g'$
are odd with respect to both $x=0$ and $x=\pi$, they vanish at such points, and $Z$ extends to a
smooth function on $[0,\pi]$, with $Z'(0)=Z'(\pi)=0$. Differentiating the first equation in
\eqref{eq:scalar} (with $\beta=0$), integrating the second one, and recalling \eqref{eq:corr-g}, 
we have
\begin{equation}\label{eq:phimg}
 -\phi''+\phi\phi'=\tau\sin x,\qquad m'=-\phi\,m,\qquad -(mg')'=m\,(A-\cos x).
\end{equation}
Finally, recall that $0<A<1$, by Lemma~\ref{lem:resp}.

We first describe the properties of $\phi$ on $(0,\pi)$.

\begin{lemma}\label{lem:shape}
On $(0,\pi)$ we have
\[
 \phi>0,\qquad \phi''\le0,\qquad (\phi\cot x)'\le0 .
\]
\end{lemma}

\begin{proof}
To start with, we show that $\phi$ is positive. Indeed, $\phi(0)=\phi(\pi)=0$ and, if 
$x_0\in(0,\pi)$ is a minimum point of $\phi$ on $[0,\pi]$, then $\phi'(x_0)=0$ and \eqref{eq:phimg} 
gives $\phi''(x_0)=-\tau\sin x_0<0$, a contradiction. Then the minimum of $\phi$ is attained only at the 
endpoints, and $\phi>0$ on $(0,\pi)$.

Next, we introduce the function
\[
 W:=\phi'\sin x-\phi \cos x ,
\]
which is smooth on $[0,\pi]$ and satisfies $W(0)=W(\pi)=0$. A direct computation, based on the 
first equation in \eqref{eq:phimg} and on its derivative, gives
\begin{equation}\label{eq:W'}
 W'=(\phi+\phi'')\sin x,
\end{equation}
and
\begin{equation}\label{eq:W}
 \tau W= \phi'\left(-\phi''+\phi\phi'\right) - \phi \left(-\phi''+\phi\phi'\right)'
 =  \phi\,\phi'''-(\phi'+\phi^2)\,\phi'' .
\end{equation}
We claim that $W\ge0$. Assume by contradiction that $W$ has a negative minimum, attained at some
$x_0\in(0,\pi)$. Then $W'(x_0)=0$, that is $\phi''(x_0)=-\phi(x_0)$, and 
\eqref{eq:W} yields
\[
 \phi'''(x_0)=-\phi'(x_0)-\phi^2(x_0)+\frac{\tau W(x_0)}{\phi(x_0)} .
\]
Differentiating \eqref{eq:W'}, and recalling that $(\phi+\phi'')(x_0)=0$, we obtain
\[
 W''(x_0)=\sin x_0\cdot\bigl(\phi'+\phi'''\bigr)(x_0)
   =\sin x_0\cdot \left(-\phi^2+\frac{\tau W}{\phi}\right)(x_0)<0,
\]
a contradiction.

We are now in a position to show that $\phi$ is concave. If not, the open set 
$\{x\in(0,\pi):\phi''(x)>0\}$ is not empty; let $(x_1,x_2)$ be one of its connected components. 
By \eqref{eq:W'}, $W'>0$ on $(x_1,x_2)$, whence $W(x_2)>0$, since $W\ge0$. In particular $x_2<\pi$, 
because $W(\pi)=0$; then $\phi''(x_2)=0$ and $\phi'''(x_2)\le0$, while \eqref{eq:W} gives
$\phi(x_2)\,\phi'''(x_2)=\tau W(x_2)>0$, a contradiction.

We are left to show that $\phi\cot x$ is nonincreasing, that is $\phi'\sin x\cos x\le\phi$ on 
$(0,\pi)$. Since $\phi$ is concave and vanishes at $0$ and $\pi$, we have $x\,\phi'(x)\le\phi(x)$ 
and $-(\pi-x)\,\phi'(x)\le\phi(x)$, whence $|\phi'(x)|\min\{x,\pi-x\}\le\phi(x)$; then
$\phi'\sin x\cos x\le|\phi'|\sin x\le\phi$, because $\sin x\le\min\{x,\pi-x\}$.
\end{proof}

Next we turn to the function $g$ or, more precisely, to the quotient $Z$.

\begin{lemma}\label{lem:Z}
On $(0,\pi)$ we have
\[
 g'>0,\qquad Z>0,\qquad Z'\ge0,\qquad (mZ)'\le0 .
\]
\end{lemma}

\begin{proof}
The function $mg'$ vanishes at $0$ and $\pi$ and, by \eqref{eq:phimg}, its derivative 
$m\,(\cos x-A)$ is positive on $(0,\arccos A)$ and negative on $(\arccos A,\pi)$ (recall that
$0<A<1$). We deduce that $mg'>0$ on $(0,\pi)$, whence $g'>0$ and $Z>0$ there.

Next, we derive the second order ODE satisfied by $Z$. Writing $g'=Z\sin x$ in the last equation of \eqref{eq:phimg}, and using $m'=-\phi m$, we obtain
\[
 \sin x\,Z'+(\cos x-\phi\sin x)\,Z=\cos x-A ;
\]
differentiating, and dividing by $\sin x$, we have
\begin{equation}\label{eq:Zeq}
 Z''+(2\cot x-\phi)\,Z'-(1+\phi'+\phi\cot x)\,Z=-1 \qquad\text{on }(0,\pi).
\end{equation}
To study the sign of the functions $Z'$ and $(mZ)'=-m(\phi Z-Z')$, it is convenient to deal with their 
weighted versions:
\[
 Y_1:=m\sin^2\!x\,Z',\qquad Y_2:=\sin^2\!x\,(\phi Z-Z').
\]
Both $Y_1$ and $Y_2$ are smooth on $[0,\pi]$, and they vanish at $0$ and $\pi$. Moreover, by a 
direct computation, \eqref{eq:Zeq} and $m'=-\phi m$ yield
\begin{equation}\label{eq:Y}
 Y_1'=m\sin^2\!x\,\bigl[(1+\phi'+\phi\cot x)\,Z-1\bigr],\qquad
 Y_2'=\sin^2\!x\,\bigl[1-(1-\phi\cot x)\,Z\bigr].
\end{equation}

Assume by contradiction that $Y_1$ has a negative minimum, attained at some $x_0\in(0,\pi)$. Then
$Z'(x_0)<0$ and $Y_1'(x_0)=0$; since $Z>0$, \eqref{eq:Y} gives 
$(1+\phi'+\phi\cot x)(x_0)=1/Z(x_0)>0$ and, differentiating,
\[
 Y_1''(x_0)=m\sin^2\!x\,\Bigl[\bigl(\phi''+(\phi\cot x)'\bigr)\,Z+(1+\phi'+\phi\cot x)\,Z'\Bigr](x_0)<0,
\]
by Lemma~\ref{lem:shape}. This is a contradiction, thus $Y_1\ge0$, that is $Z'\ge0$.

Analogously, assume that $Y_2$ has a negative minimum, attained at some $x_0\in(0,\pi)$. Then
$Z'(x_0)>\phi(x_0)Z(x_0)>0$ and $Y_2'(x_0)=0$, whence $(1-\phi\cot x)(x_0)=1/Z(x_0)>0$ and
\[
 Y_2''(x_0)=\sin^2\!x\,\Bigl[(\phi\cot x)'\,Z-(1-\phi\cot x)\,Z'\Bigr](x_0)<0,
\]
again by Lemma~\ref{lem:shape}, a contradiction. Then $Y_2\ge0$, and $(mZ)'\le0$.
\end{proof}

We are now in a position to conclude.

\begin{proof}[Proof of Proposition~\ref{prop:concavity}]
We first show that
\begin{equation}\label{eq:gsin2}
 \io m\,g\sin^2\!x\,dx>0 .
\end{equation}
To this aim, for $x\in[0,\pi/2]$ we define
\[
 \theta(x):=\frac{m(\pi-x)}{m(x)+m(\pi-x)},\qquad
 \bar g(x):=\bigl(1-\theta(x)\bigr)\,g(x)+\theta(x)\,g(\pi-x).
\]
By Lemma~\ref{lem:shape}, $m'=-\phi m<0$ on $(0,\pi)$; then the quotient $m(\pi-x)/m(x)$ has
positive derivative on $(0,\pi/2)$, and the same holds for $\theta$. Moreover
\[
 \bar g'(x)=\frac{m(x)\,g'(x)-m(\pi-x)\,g'(\pi-x)}{m(x)+m(\pi-x)}
   +\theta'(x)\,\bigl(g(\pi-x)-g(x)\bigr)>0,\qquad 0<x<\frac\pi2 :
\]
indeed, the numerator of the first term is $\sin x\,\bigl[(mZ)(x)-(mZ)(\pi-x)\bigr]$, which is 
nonnegative because $mZ$ is nonincreasing, while the second term is positive since $g$ is 
increasing (Lemma~\ref{lem:Z}). On the other hand, the change of variable $x\mapsto\pi-x$ on 
$[\pi/2,\pi]$ gives
\begin{equation}\label{eq:reflection}
 \io m\,g\,f\,dx=\int_0^{\pi/2}\bigl(m(x)+m(\pi-x)\bigr)\,\bar g(x)\,f(x)\,dx ,
\end{equation}
for every continuous function $f$ such that $f(\pi-x)=f(x)$. Since $g$ is even, 
\eqref{eq:corr-g} yields $\io mg\,dx=0$; then, taking $f\equiv1$ in \eqref{eq:reflection}, we 
infer that the strictly increasing function $\bar g$ vanishes at exactly one point $\bar x\in(0,\pi/2)$. 
Taking now $f=\sin^2\!x-\sin^2\!\bar x$, the integrand in the right hand side of 
\eqref{eq:reflection} is nonnegative, and it vanishes only at $x=\bar x$. We deduce that
\[
 \io m\,g\sin^2\!x\,dx=\io m\,g\,\bigl(\sin^2\!x-\sin^2\!\bar x\bigr)\,dx>0,
\]
that is \eqref{eq:gsin2}.

Finally, by Lemma~\ref{lem:Z}, both $g$ and $Z^2$ are nondecreasing on $[0,\pi]$, thus
\[
 \bigl(g(x)-g(y)\bigr)\bigl(Z^2(x)-Z^2(y)\bigr)\ge0\qquad\text{for every }x,y\in[0,\pi].
\]
Multiplying by $m(x)\sin^2\!x\;m(y)\sin^2\!y$, integrating on $[0,\pi]^2$, and recalling that
$Z\sin x=g'$, we obtain
\[
 \io m\sin^2\!x\,dx\cdot\io m\,g\,|g'|^2\,dx\ \ge\ \io m\,g\sin^2\!x\,dx\cdot\io m\,|g'|^2\,dx>0,
\]
where the last inequality follows from \eqref{eq:gsin2}. Since $m\,g\,|g'|^2$ is even, 
\eqref{eq:A2} yields $A''(\tau)=-6\io m\,g\,|g'|^2\,dx<0$.
\end{proof}

\end{document}